\documentclass[journal]{IEEEtran}

\usepackage{xspace,amssymb,epsfig,subfigure,syntonly}
\usepackage{epsfig,amsmath,color}
\usepackage{xspace,syntonly}
\usepackage{algorithm,algorithmic}

\newcommand{\utwi}[1]{\mbox{\boldmath $#1$}}

\newcommand{\trace}{{\textrm{Tr}}}
\newcommand{\rank}{{\textrm{rank}}}

\newcommand{\cL}{{\cal{L}}}
\newcommand{\cN}{{\cal N}}
\newcommand{\cP}{{\cal P}}

\newcommand{\cG}{{\cal G}}
\newcommand{\cA}{{\cal A}}
\newcommand{\cS}{{\cal S}}

\newcommand{\cT}{{\cal T}}

\newcommand{\cE}{{\cal E}}
\newcommand{\cI}{{\cal I}}
\newcommand{\cF}{{\cal F}}
\newcommand{\cB}{{\cal B}}

\newcommand{\cH}{{\cal H}}
\newcommand{\cV}{{\cal V}}

\newcommand{\be}{{\bf e}}

\newcommand{\br}{{\bf r}}

\newcommand{\bu}{{\bf u}}
\newcommand{\bv}{{\bf v}}
\newcommand{\bi}{{\bf i}}

\newcommand{\bA}{{\bf A}}
\newcommand{\bB}{{\bf B}}

\newcommand{\bT}{{\bf T}}
\newcommand{\bI}{{\bf I}}
\newcommand{\bX}{{\bf X}}
\newcommand{\bZ}{{\bf Z}}
\newcommand{\bW}{{\bf W}}
\newcommand{\bY}{{\bf Y}}
\newcommand{\bV}{{\bf V}}

\newcommand{\bPhi}{{\utwi{\Phi}}}

\newcommand{\bGamma}{{\utwi{\Gamma}}}
\newcommand{\bLambda}{{\utwi{\Lambda}}}

\begin{document}

\newtheorem{definition}{Definition}
\newtheorem{remark}{Remark}
\newtheorem{proposition}{Proposition}
\newtheorem{lemma}{Lemma}

%------------------------------------------------------------------------------
% Title.
%------------------------------------------------------------------------------
\title{Distributed Optimal Power Flow \\ for Smart Microgrids }

\author{Emiliano Dall'Anese, \emph{Member, IEEE},  Hao Zhu, \emph{Member, IEEE}, and Georgios B. Giannakis, \emph{Fellow, IEEE}
\thanks{\protect\rule{0pt}{0.5cm}%
Manuscript submitted on October 1, 2012; revised January 12, 2013; accepted February 16, 2013. 
This work was supported by the Inst. of Renewable Energy and the
Environment (IREE) grant no. RL-0010-13, Univ. of Minnesota. E. Dall'Anese and G. B. Giannakis are with the Digital Technology Center and the Dept. of ECE, University of Minnesota, 200 Union Street SE, Minneapolis, MN 55455, USA. Hao Zhu is with the Information Trust Institute at the University of Illinois at Urbana-Champaign, 1308 West Main St, Urbana, IL 61801, USA. 
E-mails: {\tt \{emiliano, georgios\}@umn.edu, haozhu@illinois.edu}
%
%Part of this work will be presented at the 44th North
%America Power Symposium, University of Illinois at Urbana-Champaign, IL, USA, Sep., 2012.
%
}
}

% The paper headers
\markboth{IEEE TRANSACTIONS ON SMART GRID}%
{Dall'Anese \MakeLowercase{\textit{et al.}}: Distributed Optimal Power Flow for Smart Microgrids }

% make the title area
\maketitle

%\thispagestyle{empty}

%%%%%%%%%%%%%%%%%%%%%%%%%%%%%%%%%%%%%%%%%%%%%%%
% Abstract
%%%%%%%%%%%%%%%%%%%%%%%%%%%%%%%%%%%%%%%%%%%%%%%
\begin{abstract}
Optimal power flow (OPF) is considered for microgrids, with the objective of minimizing either the power distribution losses, or, the cost of power drawn from the substation and supplied by distributed generation (DG) units, while effecting voltage regulation. The microgrid is unbalanced, due to unequal loads in each phase and non-equilateral conductor spacings on the distribution lines. Similar to OPF formulations for balanced systems, the considered OPF problem is nonconvex. Nevertheless, a semidefinite programming (SDP) relaxation technique is advocated to obtain a convex problem solvable in polynomial-time complexity. Enticingly, numerical tests demonstrate the ability of the proposed method to attain the \emph{globally} optimal solution of the original nonconvex OPF.
To ensure scalability with respect to the number of nodes, robustness to isolated communication outages, and data privacy and integrity, the proposed SDP is solved in a distributed fashion by resorting to the alternating direction method of multipliers. The resulting algorithm entails iterative message-passing among groups of consumers and guarantees faster convergence compared to competing alternatives.
\end{abstract}

\begin{keywords}
Microgrids, distribution feeders, optimal power flow, semidefinite relaxation, distributed optimization.
\end{keywords}

%%%%%%%%%%%%%%%%%%%%%%%%%%%%%%%%%%%%%%%%%%%%%%%
\section{Introduction}
\label{sec:Introduction}
%%%%%%%%%%%%%%%%%%%%%%%%%%%%%%%%%%%%%%%%%%%%%%%

Microgrids are portions of a power distribution network
located downstream of the distribution substation that supply a number of industrial and residential loads, and may include distributed generation (DG) and energy storage devices~\cite{Hatziargyriou-PESMag}. A microgrid can operate in either grid-connected, islanded, or hybrid modes. Deployment of microgrids promises drastic performance enhancement of the distribution grid in terms of efficiency and stability, along with increased network scalability and resilience to outages. 

Besides bringing power generation closer to the end user, DG units offer environment-friendly advantages over conventional generation~\cite{Hatziargyriou-PESMag}, may provide ancillary services such as reactive and harmonic compensation~\cite{Tenti08,Bolognani11}, and enable 
DG owners to actively participate in grid operations through supply contracts and pricing schemes. On the other hand, their operation must be carefully controlled in order to prevent abrupt voltage fluctuations, which stem from the well known sensitivity of voltages to variations of power injections, node over- and under-voltages~\cite{Carvalho08}, and drops of the power factor at the distribution substation~\cite{Forner}. To this end, optimal power flow (OPF) approaches are increasingly advocated also in this context to ensure efficient operation of smart microgrids and effect strict voltage regulation~\cite{Bolognani11,Forner,Tse12,Sortomme09}. 

OPF problems are deemed challenging because they require solving nonconvex problems. Nonconvexity stems from the nonlinear relationship between voltages and the complex powers demanded or injected at the nodes. In the context of transmission networks, the Newton-Raphson method~\cite{Irving87} is traditionally employed to obtain a possibly suboptimal solution of these nonconvex problems. 
In distribution systems however, its convergence is challenged by the high resistance-to-reactance ratio of distribution lines. Alternative approaches include sequential quadratic optimization, steepest descent-based methods~\cite{Forner}, fuzzy dynamic programming~\cite{Lu97}, and particle swarm optimization~\cite{Sortomme09}. However, these methods generally return suboptimal load flow solutions, and may be computationally cumbersome~\cite{PaudyalyISGT}. To alleviate these concerns, a relaxed semidefinite programming (SDP) reformulation of the OPF problem for balanced transmission systems was proposed in~\cite{Bai08} and~\cite{LavaeiLow}, where global optimality can be assessed by checking the rank of the obtained voltage-related matrix. The relaxed SDP approach was extended to \emph{balanced} distribution systems in~\cite{LavaeiGeom} and~\cite{Tse12}. Notably, for networks with a tree topology,~\cite{Tse12} and~\cite{LavaeiGeom} established sufficient conditions under which a globally optimal solution is attainable provided the original OPF problem is feasible.

Distribution networks are inherently \emph{unbalanced} because: \emph{i)} unequal single-phase loads must be served, and \emph{ii)} non-equilateral conductor spacings of three-phase line segments are involved~\cite{Kerstingbook}. 
Further, single-phase DG units may worsen the network imbalance. As a consequence, optimization approaches can not rely on single-phase equivalent models as in e.g.,~\cite{Bolognani11,Forner,Tse12,Sortomme09,LavaeiGeom}. 
For the unbalanced setup, an OPF framework was proposed in~\cite{Paudyal11}, where commercial solvers of nonlinear programs were adopted, and in~\cite{Bruno11}, where quasi-Newton methods were utilized in conjunction with load flow solvers.
However, since these methods are inherently related to gradient descent solvers of nonconvex programs, they inherit the limitations of being sensitive to initialization, and do not guarantee global optimality of their solutions. %Further, the complexity of commercial solvers of nonlinear programs do not scale well with the problem size~\cite{PaudyalyISGT}. 

The first contribution of the present paper consists in permeating the benefits of SDP relaxation techniques~\cite{luospmag10} to OPF problems for microgrids operating in an \emph{unbalanced} setup. This optimization tool has three main advantages: \emph{i)} it offers the potential of finding the \emph{globally} optimal solution; \emph{ii)} its worst-case computational complexity is quantifiable; and, \emph{iii)} it can accommodate additional thermal and quality-of-power constraints without exacerbating the problem complexity~\cite{luospmag10}. Global optimality not only reduces power distribution losses, but also leads to higher monetary savings compared to suboptimal OPF solutions.
 
The OPF problem is solved by a microgrid energy manager (MEM), which cooperates with local controllers (LCs) located throughout the network. Microgrids can vary in scope, size, and ownership~\cite{Guerrero12}. For those of medium and large size (of a distribution feeder), solving the OPF problem \emph{centrally} at the MEM may become computationally prohibitive~\cite{PaudyalyISGT}. In fact, interior point SDP solvers do not generally scale well with the problem size~\cite{Tse12,Lam11,luospmag10}. For a real-time network management, it is generally required to find a new network operational setup rapidly (e.g., in a few seconds or minutes) in order to promptly respond to abrupt load variations and to cope with the intermittent power generation that is typical of renewable-based DG units. It is then of paramount importance to solve the SDP-based OPF problem in a \emph{distributed} manner, by decomposing the main problem into multiple sub-instances that can be solved efficiently and in parallel. A distributed algorithm is desirable also when (a group of) customers do not share data with the MEM due to privacy concerns, or because they wish to manage autonomously their DG units in order to pursue individual economic interests~\cite{Guerrero12}.  Finally, a distributed algorithms involves a modest communication overhead compared to its centralized counterpart, as it does not require to pool line, generator, and load data at the MEM, and subsequently disseminate the OPF solution.  

Decentralized OPF approaches were first proposed in~\cite{Baldick97,Baldick99}, where multi-utility transmission systems were partitioned in autonomously managed areas. Augmented Lagrangian methods were employed to decompose the overall OPF problem in per-area instances. A similar approach was followed by e.g.,~\cite{Hug09,Nogales03} (see also references therein), where standard Lagrangian approaches were utilized in conjunction with Newton methods. Solving SDP in a distributed fashion is challenging due to the couplings of local voltage-related matrices enforced by the positive semidefinite constraint of the global voltage matrix. Results related to positive semidefinite matrix completion~\cite{Grone84} were leveraged in~\cite{Tse12,Lam11} to develop a distributed OPF algorithms for balanced networks via dual decomposition. Using the results of~\cite{Grone84}, and tapping into the powerful alternating direction method of multipliers (ADMM)~\cite[Sec.~3.4]{BeT89}, a distributed optimization problem for unbalanced microgrids is formulated here, where each LC solves an optimization sub-problem, and then exchanges simple messages with its neighboring LCs. Compared to the dual decomposition schemes of~\cite{Tse12,Lam11}, the proposed ADMM-based approach offers a markedly improved convergence. 

In the OPF context, augmented Lagrangian methods (related to ADMM) were first used in~\cite{Baldick97,Baldick99} to develop a decentralized optimization scheme for (balanced) transmission networks. Off-the-shelf schemes were used to solve the OPF sub-problem associated with each sub-network. More recently, ADMM was advocated in~\cite{Kekatos13} for state estimation, and in~\cite{Boyd_messPass} for distributed multi-period OPF in balanced systems. Here, the approaches of~\cite{Baldick97,Baldick99,Boyd_messPass} are considerably broadened by considering unbalanced distribution networks, and by employing the ADMM to devise a distributed SDP solver.    

The rest of the paper is organized as follows. Section~\ref{sec:Modeling}
recapitulates the OPF formulation for microgrids, and Section
\ref{sec:semidefinite} develops its centralized SDP solver. 
The distributed algorithm is presented in Section~\ref{sec:distributed}, while numerical tests are reported in Section~\ref{sec:results}. Finally, concluding remarks are provided in Section~\ref{sec:conclusions}\footnote{{\it Notation:} Upper (lower) boldface
letters will be used for matrices (column vectors); $(\cdot)^\cT$ for transposition; $(\cdot)^*$ complex-conjugate; and, $(\cdot)^\cH$ complex-conjugate transposition;  
$\Re\{\cdot\}$ denotes the real part, and $\Im\{\cdot\}$ the imaginary part; $j := \sqrt{-1}$ the imaginary unit. $\trace(\cdot)$ the matrix trace; $\rank(\cdot)$ the matrix rank; $|\cdot|$ denotes the magnitude of a number or the cardinality of a set; and $\|\cdot\|_F$ stands for Frobenius norm.  
Given a vector $\bv$ and a matrix $\bV$, $[\bv]_{\cP}$ denotes a $|\cP| \times 1$ sub-vector containing entries of $\bv$ indexed by the set $\cP$, and $[\bV]_{\cP_1,\cP_2}$ the $|\cP_1| \times |\cP_2|$ sub-matrix with row and column indexes described by $\cP_1$ and $\cP_2$. Finally, $\mathbf{0}_{M\times N}$ and $\mathbf{1}_{M\times N}$ denotes $M \times N$ matrices with all zeroes and ones, respectively.}.

%%%%%%%%%%%%%%%%%%%%%%%%%%%%%%%%%%%%%%%%%%%%%%%
\section{Modeling and problem formulation}
\label{sec:Modeling}
%%%%%%%%%%%%%%%%%%%%%%%%%%%%%%%%%%%%%%%%%%%%%%%

Consider a microgrid comprising $N$ nodes collected
in the set\footnote{The symbols defined throughout the paper are recapitulated in Table~\ref{tab:nomenclature}.} $\cN := \{0, 1,\ldots,N\}$, and overhead or underground lines represented by the set of edges $\cE := \{(m,n)\} \subset \cN  \times \cN$. Let node $0$ represent the point of common coupling (PCC), taken to be the distribution substation. Define as $\cP_{mn} \subseteq \{a_{mn},b_{mn},c_{mn}\}$ and $\cP_{n} \subseteq \{a_{n},b_{n},c_{n}\}$ the phases of line $(m,n) \in \cE$ and node $n \in \cN$, respectively. Let $V_n^{\phi} \in \mathbb{C}$ be the complex line-to-ground voltage at node $n \in \cN$ of phase $\phi \in \cP_n$, and $I_n^{\phi} \in \mathbb{C}$ the current injected at the same node and phase. As usual, the voltages $\bv_0 :=  [V_0^{a},V_0^{b},V_0^{c}]^{\cT}$ are taken as reference for the phasorial representation. 
Lines $(m,n) \in \cE$ are modeled as $\pi$-equivalent components~\cite[Ch. 6]{Kerstingbook}, and the $|\cP_{mn}| \times |\cP_{mn}|$ phase impedance and shunt admittance matrices are denoted as $\bZ_{mn} \in \mathbb{C}^{|\cP_{mn}| \times |\cP_{mn}|}$ and $\bY_{mn}^{(s)} \in \mathbb{C}^{|\cP_{mn}| \times |\cP_{mn}|}$, respectively. 
%If four-wire grounded wye lines or lines with multi-grounded neutrals are present, the $|\cP_{mn}| \times |\cP_{mn}|$ matrices $\bZ_{mn}$ and $\bY_{mn}^{(s)}$ can be obtained from the higher-dimensional ``primitive'' matrices via Kron reduction~\cite[Ch. 6]{Kerstingbook}. 
Three- or single-phase transformers (if any) are modeled as series components with transmission parameters that depend on the connection type~\cite[Ch.~8]{Kerstingbook},~\cite{Paudyal11}.

Per phase $\phi \in \cP_n$, let $P_{L,n}^{\phi}$ and $Q_{L,n}^{\phi}$ denote the active and reactive powers demanded by a wye-connected load at the bus $n$. %In case of delta-connected load, a delta-to-wye transformation can be typically employed. 
Capacitor banks are usually mounted at some nodes to provide reactive power support. Let $y_{C,n}^{\phi}$ denote the susceptance of a capacitor connected at node $n$ and phase $\phi$. 
Finally, suppose that $S$ DG units are located at a subset of nodes $\cS \subset \cN$, and let $P_{G,s}^{\phi}$ and $Q_{G,s}^{\phi}$ denote the active and reactive powers supplied by unit $s \in  \cS$. For conventional DG units, such as diesel generators, the supplied powers can be controlled, and they will be variables of the OPF problem; on the other hand, $P_{G,s}^{\phi}$ and $Q_{G,s}^{\phi}$ represent committed powers for renewable-based sources. 

Given the demanded loads $\{P_{L,n}^{\phi},Q_{L,n}^{\phi}\}$, 
the goal is to select a feasible set of voltages $\{V_{n}^{\phi}\}$, currents $\{I_{n}^{\phi}\}$, and powers supplied by conventional DG units $\{P_{G,s}^{\phi},Q_{G,s}^{\phi}\}$ so that the steady-state operation of the microgrid is optimal in a well defined sense. To this end, one of the following two objectives is usually pursued.

\noindent \emph{i) Minimization of power losses}. The active power flowing on line $(m,n) \in \cE$ is $P_{m\rightarrow n} := \sum_{\phi \in \cP_{mn}} V_{m}^\phi (I_{m\rightarrow n}^{\phi})^*$, where $I_{m\rightarrow n}^{\phi}$ is the complex current from node $m$ to node $n$ on phase $\phi \in \cP_{mn}$. Thus, the cost to be minimized here is~\cite{Tse12}
\begin{align}
C_1(\cV) = \sum_{(m,n) \in \cE} \left(  P_{m\rightarrow n} + P_{n\rightarrow m} \right)
\label{power_loss_cost}  
\end{align}
where $\cV := \{I_0^{\phi}, \{V_{n}^{\phi}\}, \{I_n^{\phi}\}, \{P_{G,n}^{\phi},Q_{G,n}^{\phi}\}\}$ collects all the steady-state variables.

\noindent \emph{ii) Minimization of the cost of supplied power}. Letting the cost of power drawn from the PCC be denoted by $c_0 > 0$, and the one incurred by the use of DG unit $s \in \cS$ as $c_s \geq 0$, one can minimize the cost of supplied power~\cite{Forner}
\begin{equation}
C_2(\cV) :=  c_0 \sum_{\phi \in \cP_{0}} V_{0}^\phi (I_{0}^\phi)^* + \sum_{s \in \cS} c_{s} \sum_{\phi \in \cP_{s}} P_{G,s}^{\phi} \, .
\label{econ_dispatch_cost}  
\end{equation}
Notice that~\eqref{power_loss_cost} and~\eqref{econ_dispatch_cost} are equivalent when $c_0 = 1$ and $c_s = 1$ for all $s \in \cS$.

Based on~\eqref{power_loss_cost}--\eqref{econ_dispatch_cost}, the following OPF problem is considered:
\vspace{-.3cm}
\begin{subequations} 
\label{Pmg}
\begin{align} 
& \hspace{-.8cm} \mathrm{(P1)} \hspace{1cm} \min_{\cV} \,\, C_m(\cV) \label{mg-cost} \\
\textrm{s.t.} \,\, 
& V_s^{\phi} (I_s^{\phi})^* = P_{G,s}^{\phi} - P_{L,s}^{\phi} + j (Q_{G,s}^{\phi} - Q_{L,s}^{\phi})  , \nonumber \\
& \hspace{4.2cm}  \forall \,\,\phi \in \cP_s, \, s \in \cS \label{mg-balance-source} \\
& V_n^{\phi} (I_n^{\phi})^* = - P_{L,n}^{\phi} - jQ_{L,n}^{\phi} + j y_{C,n}^{\phi} |V_n^{\phi}|^2,  \nonumber \\ 
& \hspace{4.2cm}  \forall\,\, \phi \in \cP_n, \,  n \in \cN \backslash \cS \label{mg-balance} \\
& I_n^{\phi}  = \sum_{m \in \cN_m}  \Big[ \Big(\frac{1}{2}\bY_{mn}^{(s)} + \bZ_{mn}^{-1} \Big) [\bv_n]_{\cP_{mn}}  \nonumber  \\
& \hspace{1.65cm} - \bZ_{mn}^{-1} [\bv_m]_{\cP_{mn}} \Big]_{\{\phi\}} \forall \,\, \phi \in \cP_n,\,  n \in \cN  \label{mg-current} \\
& V_n^{\mathrm{min}} \leq |V_n^{\phi}| \leq V_n^{\mathrm{max}} , \hspace{1.2cm} \forall \,\, \phi \in \cP_n,\,  n \in \cN  \label{mg-Vlimits} \\
& P_{G,s}^{\textrm{min}} \leq P_{G,s}^{\phi} \leq  P_{G,s}^{\textrm{max}}, \hspace{1.4cm} \forall \,\, \phi \in \cP_s,\,  s \in \cS \label{mg-plimits} \\
& Q_{G,s}^{\textrm{min}} \leq Q_{G,s}^{\phi} \leq  Q_{G,s}^{\textrm{max}}, \hspace{1.3cm} \forall \,\, \phi \in \cP_s,\,  s \in \cS
\label{mg-qlimits} 
\end{align}
\end{subequations}
where $V_n^{\mathrm{min}}$ and $V_n^{\mathrm{max}}$ in~\eqref{mg-Vlimits} are given minimum and maximum utilization and service voltages; $\cN_n := \{j|(n,j) \in \cE\}$;~\eqref{mg-plimits}--\eqref{mg-qlimits} are box constraints for the power supplied by conventional DG units, and $m \in \{1,2\}$ depending on the chosen cost.

Similar to OPF variants for transmission networks and balanced distribution networks, (P1) is a nonlinear \emph{nonconvex} problem because of the load flow equations~\eqref{mg-balance-source}-\eqref{mg-balance} as well as the constraints~\eqref{mg-Vlimits}. In the next section, an equivalent reformulation of (P1) will be introduced, and its solution will be pursued using the SDP relaxation technique.

%%%%%%%%%%%%%%%%%%%%%%%%
\begin{table}[t]
\caption{Nomenclature}
\begin{center}
\begin{tabular}{ll}
 \hline  \\ 
$\cN$ & Set collecting the nodes of the microgrid \\
$\cE$ & Set collecting the distribution lines of the microgrid \\
$\cS$ & Subset of nodes featuring DG  \\
$\cN_n$ & Set of neighboring nodes of $n \in \cN$ \\
$\cP_n$ & Set of phases at node $n$  \\
$\cP_{mn}$ & Set of phases of line $(m,n)$  \\
$V_n^\phi$ & Complex line-to-ground voltage at phase $\phi$ of node $n$   \\
$I_n^\phi$ & Complex current injected at phase $\phi$ of node $n$   \\
$I_{m \rightarrow n}^\phi$ & Complex current on phase $\phi$ of line $(m,n)$    \\
$P_{L,n}^{\phi} (Q_{L,n}^{\phi})$ & Active (reactive) load demanded at node $n$ on phase $\phi$  \\
$P_{G,n}^{\phi} (Q_{G,n}^{\phi})$ & Active (reactive) power supplied at node $n$ on phase $\phi$  \\
$y_{C,n}^{\phi}$ & Susceptance of a capacitor at node $n$ and phase $\phi$ \\
$P_{m \rightarrow n}$ & Active power exiting node $m$ on line $(m,n)$ \\
$\bZ_{mn}$ & Phase impedance matrix of line $(m,n)$ \\
$\bY^{(s)}_{mn}$ & Shunt admittance matrix of line $(m,n)$ \\
$V_n^{\mathrm{min}}$, $V_n^{\mathrm{max}}$ & Utilization and service voltage magnitude limits \\
$P_{G,n}^{\textrm{min}} (Q_{G,n}^{\textrm{min}})$ & Minimum active (reactive) power supplied at node $n$ \\
$P_{G,n}^{\textrm{max}} (Q_{G,n}^{\textrm{max}})$ & Maximum active (reactive) power supplied at node $n$ \\
$\cA^{(\ell)}$ & Set of nodes forming the area $\ell$ of the microgrid \\
$\bar{\cA}^{(\ell)}$ & Extended area $\ell$ comprising the nodes in $\cA^{(\ell)}$ and the \\
& nodes of different areas connected to $\cA^{(\ell)}$ by a line \\
$\bar{\cN}^{(\ell)}$ & Set of neighboring areas of $\cA^{(\ell)}$ \\
\hline 
\end{tabular}
\end{center}
\label{tab:nomenclature}
\end{table}%
%%%%%%%%%%%%%%%%%%%%%%%%

%%%%%%%%%%%%%%%%%%%%%%%%%%%%%%%%%%%%%%%%%%%%%%%
\section{SDP-Based Centralized Solution}
\label{sec:semidefinite}
%%%%%%%%%%%%%%%%%%%%%%%%%%%%%%%%%%%%%%%%%%%%%%%

Consider the $\sum_{n\in\cN} |\cP_n| \times 1$ complex vectors $\bv := [\bv_0^\cT,\bv_1^\cT, \ldots,\bv_{N}^\cT]^\cT$ and
$\bi := [\bi_0^\cT,\bi_1^\cT, \ldots,\bi_{N}^\cT]^\cT$, with $\bv_n$ and $\bi_n$ the $|\cP_n| \times 1$ complex vectors collecting voltages $\{V_{n}^\phi\}_{\phi \in \cP_n}$ and currents $\{I_{n}^\phi\}_{\phi \in \cP_n}$ per node $n \in \cN$. 
Voltages and injected currents abide by Ohm's law $\bi = \bY \bv$, where $\mathbf{Y}$ is a symmetric block matrix of dimensions $\sum_{n\in\cN} |\cP_n| \times \sum_{n\in\cN} |\cP_n|$, whose entries are given by: 

\emph{i)} matrix $-\bZ_{mn}^{-1}$ occupying the $|\cP_{mn}| \times |\cP_{mn}|$ off-diagonal block corresponding to line $(m,n) \in \cE$; and, 

\emph{ii)} the $|\cP_{n}| \times |\cP_{n}|$ diagonal block corresponding to node $n \in \cN$ with $\cN_m := \{n|(m,n) \in \cE\}$

\begin{align} 
[\bY]_{\cP_n,\cP_n} := \sum_{m \in \cN_n}  \left(\frac{1}{2}\tilde{\bY}_{mn}^{(s)} + \tilde{\bY}_{mn} \right) \label{Ymatrix}
\end{align}
where $\tilde{\bY}_{mn} := \bZ_{mn}^{-1}$ if $\cP_n = \cP_{mn}$, otherwise $[\tilde{\bY}_{mn}]_{\cP_{nm},\cP_{nm}} := \bZ_{mn}^{-1}$ and $[\tilde{\bY}_{mn}]_{\cP_n \backslash \cP_{nm},\cP_n \backslash \cP_{nm}} = \mathbf{0}$ ($\tilde{\bY}_{mn}^{(s)}$ is formed likewise). 

The next step consists in expressing the active and reactive powers injected per node, active powers flowing on the lines, and voltage magnitudes, as \emph{linear} functions of the outer-product matrix $\bV := \bv \bv^\cH$. To this end, define the following admittance-related matrix per node $n$ and phase $\phi$
\begin{equation}
\bY_n^{\phi} := \bar{\be}_n^{\phi} (\bar{\be}_n^{\phi})^{\cT} \bY 
\label{eq:Ynode}
\end{equation}
where $\bar{\be}_n^{\phi} := [\mathbf{0}_{|\cP_0|}^{\cT},\ldots,\mathbf{0}_{|\cP_{n-1}|}^{\cT},\be_{\cP_n}^{\phi, \cT},\mathbf{0}_{|\cP_{n+1}|}^{\cT},\ldots,\mathbf{0}_{|\cP_N|}^{\cT}]^{\cT}$, and $\{\be_{\cP_n}^\phi\}_{\phi \in \cP_n}$ denotes the canonical basis of $\mathbb{R}^{|\cP_n|}$, and let the $|\cP_{mn}| \times \sum_{n\in\cN} |\cP_n|$ matrices $\bA_{m \rightarrow n}$ and $\bB_{m}$ be defined as 
$\bA_{m \rightarrow n} := [\mathbf{0}_{|\cP_{mn}|\times \sum_{i = 0}^{m-1}|\cP_i|},$ $\bZ_{mn}^{-1},\mathbf{0}_{|\cP_{mn}|\times\sum_{i = m+1}^{n-1}|\cP_i|},-\bZ_{mn}^{-1},$ $\mathbf{0}_{|\cP_{mn}|\times\sum_{i = n+1}^{N}|\cP_i|}]$ and $\bB_{m} := [\mathbf{0}_{|\cP_{mn}|\times \sum_{i = 0}^{m-1}|\cP_i|},\bI_{|\cP_{mn}|}, \mathbf{0}_{|\cP_{mn}|\times\sum_{i = m+1}^{N}|\cP_i|}]$, respectively. 
Then, a linear model in $\bV$ can be established using the following lemma (see also~\cite{ZhuNAPS11} and~\cite{LavaeiLow}).

\vspace{.2cm}
\begin{lemma}
\label{SDPreformulation}
For the Hermitian matrices   
\begin{subequations}
\label{eq:Phi}
\begin{align}
\bPhi_{P,n}^{\phi} & := \frac{1}{2} (\bY_n^{\phi} + (\bY_n^{\phi})^\cH) \label{eq:PhiP} \\
\bPhi_{Q,n}^{\phi} & := \frac{j}{2} (\bY_n^{\phi} - (\bY_n^{\phi})^\cH)  \label{eq:PhiQ} \\
\bPhi_{V,n}^{\phi} & := \bar{\be}_n^{\phi} (\bar{\be}_n^{\phi})^{\cT}  \label{eq:PhiV} \\
\bPhi_{m\rightarrow n} & :=  \frac{1}{2}\left(\bA_{m\rightarrow n}^\cH \bB_{m}  + \bB_{m}^\cH \bA_{m \rightarrow n}\right) \label{eq:PhiPmn}
\end{align}
\end{subequations}
voltage magnitudes and active as well as reactive powers are linearly related with $\bV$ as
\begin{subequations}
\label{eq:PQV} 
\begin{align}
\trace(\bPhi_{V,n}^{\phi} \bV) & = |V_{n}^{\phi}|^2  \label{eq:V} \\
\trace(\bPhi_{P,n}^{\phi} \bV) & = P_{G,n}^{\phi} - P_{L,n}^{\phi} \label{eq:P} \\
\trace(\bPhi_{Q,n}^{\phi} \bV) & = Q_{G,n}^{\phi} - Q_{L,n}^{\phi} + y_{C,n}^{\phi}\trace(\bPhi_{V,n}^{\phi} \bV)  \label{eq:Q}  \\
\trace(\bPhi_{m\rightarrow n} \bV) & = P_{m\rightarrow n}  \label{eq:Pmn}
\end{align}
\end{subequations}
with $P_{G,n}^{\phi} = Q_{G,n}^{\phi} = 0$ for $n \in \cN\backslash \cS$, and $y_{C,n}^{\phi} = 0$ if capacitor banks are not present at node $n$. 
\end{lemma}
\emph{Proof.}   See the Appendix.  \hfill $\Box$

\vspace{.2cm}

Using Lemma~\ref{SDPreformulation}, problem (P1) is \emph{equivalently} reformulated as follows: 
\begin{subequations}
\label{mg2}
\begin{align} 
& \hspace{-.7cm} \mathrm{(P2)} \hspace{1cm} \min_{\bV} \,\, \tilde{C}_m(\bV) \nonumber  \\
\textrm{s.t.} \,\, 
& \trace(\bPhi_{P,n}^{\phi} \bV) + P_{L,n}^{\phi} = 0, \hspace{1cm} \forall\,\, \phi, \, \forall \, n \in \cN \backslash \cS  \label{mg2-P} \\
& \trace(\bPhi_{Q,n}^{\phi} \bV) + Q_{L,n}^{\phi} - y_{C,n}^{\phi}\trace(\bPhi_{V,n}^{\phi} \bV)  = 0, \nonumber \\
& \hspace{4.6cm} \forall\,\, \phi, \, \forall \, n \in \cN \backslash \cS \label{mg2-Q} \\
&P_{G,s}^{\textrm{min}} \leq \trace(\bPhi_{P,s}^{\phi} \bV) + P_{L,s}^{\phi} \leq P_{G,s}^{\textrm{max}},  \,\,\, \forall\,\, \phi, \, \forall \, n \in \cS  \label{mg2-Pg} \\
&Q_{G,s}^{\textrm{min}} \leq \trace(\bPhi_{Q,s}^{\phi} \bV) + Q_{L,s}^{\phi} \leq Q_{G,s}^{\textrm{max}},  \, \forall\,\, \phi, \, \forall \, n \in \cS  \label{mg2-Pg} \\
& (V_n^{\mathrm{min}})^2 \leq \trace(\bPhi_{V,n}^{\phi} \bV) \leq (V_n^{\mathrm{max}})^2, \,\, \,\, \forall\,\, \phi, \, \forall \, n \in \cN \label{mg2-voltage} \\
& \bV \succeq \mathbf{0} , \,\, [\bV]_{\cP_0,\cP_0} = \bv_0 \bv_0^{\cH} \label{mg2-semipos}  \\
& \rank(\bV) = 1   \label{mg2-rank} 
\end{align}
\end{subequations}
where the costs $\tilde{C}_1(\bV)$ and $\tilde{C}_2(\bV)$ are re-expressed as
\begin{subequations}
\begin{align} 
\tilde{C}_1(\bV) &  = \sum_{(m,n) \in \cE} \left( \trace(\bPhi_{m \rightarrow n} \bV) + \trace(\bPhi_{n \rightarrow m} \bV) \right)  \label{eq:refcost1} \\
\tilde{C}_2(\bV) & =  \sum_{s \in \cS \cup \{0\}}  c_s \sum_{\phi \in \cP_s} \trace(\bPhi_{P,s}^{\phi} \bV) \, . \label{eq:refcost2}  
\end{align}
\end{subequations}

Problem (P2) is still nonconvex because of the rank-1 constraint~\eqref{mg2-rank}. Nevertheless, the SDP relaxation technique, which amounts to dropping the rank constraint~\cite{luospmag10}, can be leveraged to obtain the following \emph{convex} problem: 
\begin{align} 
& \hspace{-3.5cm} \mathrm{(P3)} \hspace{2cm} \min_{\bV} \,\, \tilde{C}_m(\bV)   \\
\textrm{s.t.} & \,\, \eqref{mg2-P}-\eqref{mg2-semipos} \, . \nonumber
\end{align}
Clearly, if the optimal solution $\bV_{\textrm{opt}}$ of (P3) has rank 1, then 
it is a \emph{globally} optimal solution also for the nonconvex problem (P2). Further, since (P1) and (P2) are equivalent, there exists a vector $\bv_{\textrm{opt}}$ so that $\bV_{\textrm{opt}} = \bv_{\textrm{opt}} \bv_{\textrm{opt}}^\cH$, and the optimal costs of (P1) and (P2) coincide at the optimum. This is formally summarized next.  

\vspace{.2cm}

\begin{proposition}  
Let $\bV_{\textrm{opt}}$ denote the optimal solution of the SDP (P3), and assume that $\rank(\bV_{\textrm{opt}}) = 1$. Then, the vector of line-to-ground voltages
%\begin{align} 
$\bv_{\textrm{opt}} := \sqrt{\lambda_1} \bu_1$,
%\end{align}
where $\lambda_1 \in \mathbb{R}^+$ is the unique non-zero eigenvalue of $\bV_{\textrm{opt}}$ and $\bu_1$ the corresponding eigenvector, is a globally optimal solution of (P1).  
\hfill $\Box$
\end{proposition}

\vspace{.2cm}

The upshot of the proposed formulation is that a \emph{globally} optimal solution of (P2) (and hence (P1)) can be obtained via standard interior-point solvers in \emph{polynomial time}. In fact, the worst-case complexity of (P3) is on the order $\mathcal{O}(\max\{N_c,\sum_{n\in\cN} |\cP_n|\}^4 \sqrt{\sum_{n\in\cN} |\cP_n|} \log(1/\epsilon))$ for general purpose SDP solvers, with $N_c$ denoting the total number of constraints and $\epsilon > 0$ a given solution accuracy~\cite{luospmag10}. Notice however that the sparsity of $\{\bPhi_{P,n}^{\phi},\bPhi_{Q,n}^{\phi}, \bPhi_{V,n}^{\phi}\}$ and the so-called chordal structure of a matrix can be exploited to obtain substantial computational savings; see e.g.,~\cite{Jabr12}. 
In contrast, gradient descent-based solvers for nonconvex programs, sequential quadratic programming, and particle swarm optimization, do not guarantee global optimality of the obtained solutions and are sensitive to initialization. Here, global optimality translates to lower distribution losses and increased monetary savings compared to sub-optimal OPF solutions. 

Since (P3) is a relaxed version of (P2), $\bV_{\textrm{opt}}$ could have rank greater than $1$. In this case, rank reduction techniques may be employed to find a feasible rank-1 approximation of $\bV_{\textrm{opt}}$ provided it exists. For instance, the randomization technique offers a viable way to obtain a rank-1 approximation with quantifiable approximation error; see e.g.,~\cite{luospmag10} and references therein. Albeit feasible for (P2), the resultant solution is generally suboptimum~\cite{luospmag10}. For \emph{balanced} tree distribution networks,~\cite{LavaeiGeom} and~\cite{Tse12} established conditions under which a rank-$1$ solution is attainable provided the original OPF problem is feasible. Unfortunately, when the tree power network is \emph{unbalanced}, the results of~\cite{LavaeiGeom} and~\cite{Tse12} no longer apply, as explained in the ensuing Section~\ref{sec:rank}. However, an intuitive argument will be provided in Section~\ref{sec:rank} to explain why a rank-$1$ solution is expected even in the unbalanced setup. But first, SDP-consistent constraints on line flows are derived in Section~\ref{sec:current}, and a remark is provided.

\vspace{.1cm}

\noindent \emph{Remark 1}. Step-down or in-line three- or single-phase transformer banks (if any) can be accommodated in the formulated optimization problems by using their series component models~\cite[Ch.~8]{Kerstingbook},~\cite{Paudyal11}. If a delta connection is employed on one side of the transformer, a small ``dummy'' resistance should be included between the primary and the secondary sides (one per phase) in order to ensure that the matrix $\bV_{\textrm{opt}}$ obtained by solving (P3) has rank $1$; see also~\cite{LavaeiLow}.  \hfill $\Box$

%%%%%%%%%%%%%%%%%%
\subsection{Constraints on line flows}
\label{sec:current}
%%%%%%%%%%%%%%%%%%
Constraints on the power dispelled on the distribution lines, or, on the line current magnitudes are generally adopted to protect conductors from overheating (which may eventually trigger an outage event). 
Using Lemma~\ref{SDPreformulation}, it turns out that the real power dissipated on a line $(m,n) \in \cE$ can be limited by simply adding the constraint $\trace(\bPhi_{m \rightarrow n} \bV) + \trace(\bPhi_{n \rightarrow m} \bV)  \leq \Delta P_{mn}$ in (P3), for a given maximum power loss $\Delta P_{mn}$. 

Consider now the constraint $|I_{mn}^{\phi}| \leq  I_{mn}^{\textrm{max}}$, with $I_{mn}^{\textrm{max}}$ a given upper bound on the magnitude of $I_{mn}^{\phi}$.  Aiming to an SDP-consistent reformulation of this constraint, let $\bi_{mn} : = [\{I_{mn}^{\phi}\}]^\cT$ denote the $|\cP_{mn}| \times 1$ vector collecting the complex currents flowing through line $(m,n) \in \cE$, and notice that $\bi_{mn}$ is related to voltages $\bv_{m}$ and $\bv_{n}$ as $\bi_{mn}  =  \bZ_{mn}^{-1} \left([\bv_{m}]_{\cP_{mn}} - [\bv_{n}]_{\cP_{mn}} \right)$. Next, define the $|\cP_{mn}| \times \sum_{n \in \cN} |\cP_n|$ complex matrix
\begin{align}
\bB_{mn} & := [\mathbf{0}_{|\cP_{mn}| \times \sum_{n=0}^{m-1} |\cP_n|}, \check{\bZ}_{mn}^m, \ldots \nonumber \\
& \mathbf{0}_{|\cP_{mn}| \times \sum_{n=m+1}^{n-1} |\cP_n|}, \check{\bZ}_{mn}^n  \mathbf{0}_{|\cP_{mn}| \times \sum_{n=n+1}^{N} |\cP_n|}]
\end{align}
where $\check{\bZ}_{mn}^m$ is a $|\cP_{mn}| \times |\cP_{m}|$ matrix with elements $[\check{\bZ}_{mn}^m]_{\cP_{mn},\cP_{mn}} = \bZ_{mn}^{-1}$ and $[\check{\bZ}_{mn}^m]_{\cP_{mn},\cP_m \backslash \cP_{mn}} = \mathbf{0}$; likewise, $\check{\bZ}_{mn}^n$ has dimensions $|\cP_{mn}| \times |\cP_{n}|$, and its entries are $[\check{\bZ}_{mn}^n]_{\cP_{mn},\cP_{mn}} = -\bZ_{mn}^{-1}$ and $[\check{\bZ}_{mn}^n]_{\cP_{mn},\cP_n \backslash \cP_{mn}} = \mathbf{0}$. Thus, an SDP-compliant re-formulation of the constraint on the current magnitude is possible as follows.   

\begin{lemma}
\label{limitcurrent}
Consider the Hermitian matrix 
\begin{equation}
\bPhi_{I,mn}^{\phi} :=  \bB_{mn}^{\cH} \be_{mn}^{\phi} (\be_{mn}^{\phi})^{\cT} \bB_{mn}
\end{equation}
where $\{\be_{mn}^\phi\}_{\phi \in \cP_{mn}}$ denotes the canonical basis of $\mathbb{R}^{|\cP_{mn}|}$. Then, constraint $|I_{mn}^{\phi}| \leq  I_{mn}^{\textrm{max}}$ can be equivalently re-expressed as   
\begin{align}
\trace\{ \bPhi_{I,mn}^{\phi} \bV \} & \leq  (I_{mn}^{\textrm{max}})^2 \, .  \label{currentSDP} 
\end{align}
\end{lemma}
\hfill $\Box$

Following similar steps, and using~\eqref{mg-current}, constraints on the magnitude of the injected currents $\{I_n^{\phi}\}$ can be derived too. 

In unbalanced microgrids, it is of prime interest to protect from overheating also the neutral cable(s) of the distribution lines~\cite{Kerstingbook}. Towards this end, let $\cP_{mn}^{(\varphi)}$ denote the set of grounded neutral cables that are present on the line $(m,n) \in \cE$. Further, let $\bT_{mn}$ the  $|\cP_{mn}^{(\varphi)}| \times |\cP_{mn}|$ be the neutral transformation matrix, which is obtained from the primitive impedance matrix of line $(m,n)$ via Kron reduction~\cite[Sec.~4.1]{Kerstingbook}. Thus, the neutral currents $\bi^{(\varphi)}_{mn} := [I^{(1)}_{mn}, \ldots, I^{|\cP_{mn}^{\varphi}|}_{mn}]^{\cT}$ are linearly related to the line currents $\bi_{mn}$ as $\bi^{(\varphi)}_{mn} = \bT_{mn} \bi_{mn}$. It readily follows from Lemma~\ref{limitcurrent}, that the magnitude of the current on the neutral cables can be constrained in (P3) as  
\begin{equation}
\trace\{ \bPhi_{I,mn}^{(\varphi)} \bV \}  \leq  (I_{mn}^{(\varphi), \textrm{max}})^2,    \, \quad  \forall \, \varphi \in \cP^{(\varphi)}_{mn}   \label{neutralcurrentlimit}
\end{equation}
with $\bPhi_{I,mn,t}^{(\varphi)} := \bB_{mn}^{\cH} \bT_{mn}^{\cH} \be_{mn}^{(\varphi)} (\be_{mn}^{(\varphi)})^{\cT} \bT_{mn} \bB_{mn}$, $\{\be_{mn}^{(\varphi)}\}$ the canonical basis of $\mathbb{R}^{|\cP^{(\varphi)}_{mn}|}$, and $I_{mn}^{(\varphi), \textrm{max}}$ a cap on the magnitude of $I_{mn}^{(\varphi)}$.

%%%%%%%%%%%%%%%%%%
\subsection{The rank conundrum}
\label{sec:rank}
%%%%%%%%%%%%%%%%%%

Sufficient conditions under which a rank-$1$ solution is always obtained provided the original OPF problem is feasible were established in~\cite{LavaeiGeom} and~\cite{Tse12} for \emph{balanced} distribution networks with a \emph{tree} topology. Balanceness implies that equal single-phase loads are served, $\bZ_{mn}^{-1} = (g_{mn} - j b_{mn})\bI_{|\cP_{nm}|}$ and $\bY_{mn}^{(s)} = j b^{(s)}_{mn}\bI_{|\cP_{nm}|}$ for each line $(m,n) \in \cE$, where $g_{mn}, b_{mn}, b^{(s)}_{mn} > 0$. To recapitulate the broad outline of the proofs in~\cite{LavaeiGeom} and~\cite{Tse12}, assume for simplicity that the shunt admittances are all zero. Then, the total power flowing from node $m$ to $n$ is given by $P_{m\rightarrow n} = 3|V_m^{\phi}|^2 + 3|V_m^{\phi}||V_n^{\phi}| (b_{mn} \sin(\theta_{mn}^{\phi}) -  g_{mn} \cos(\theta_{mn}^{\phi}))$, with $\theta_{mn}^{\phi} := \theta_{m}^{\phi} - \theta_{n}^{\phi}$ the angle difference between voltages $V_m^{\phi}$      
and $V_n^{\phi}$. Since the network is balanced, $\theta_{mn}^{\phi}$ is the same on each phase $\phi \in \cP_{mn}$. Fixing the voltage magnitudes, the region of the feasible powers $(P_{m\rightarrow n},P_{n\rightarrow m} )$, which is denoted as $\cF_{mn}$, becomes an affine transformation of the unit circles. Then, if one minimizes a strictly increasing function of the powers $(P_{m\rightarrow n},P_{n\rightarrow m})$, it follows that the Pareto front of $\cF_{mn}$ and the one of the convex hull of $\cF_{mn}$ coincide if $-\tan^{-1}(b_{nm}/g_{nm}) < \theta_{nm}^{\phi} < \tan^{-1}(b_{nm}/g_{nm})$~\cite{LavaeiGeom}.  
Based on this observation, proving that the Pareto regions of the feasible power injections at the nodes of the nonconvex OPF (the balanced counterpart of (P2)) and the relaxed SDP (the balanced counterpart of (P3)) amounts to showing that: \emph{i)} the region of feasible powers $\cF := \{(P_{m\rightarrow m},P_{m\rightarrow m})\}_{(m,n)\in \cE}| (P_{m\rightarrow m},P_{m\rightarrow m}) \in \cF_{mn}  \,\, \forall \,\, (m,n)\in \cE\}$ is the Cartesian product of the regions $\{\cF_{mn}\}_{(m,n)\in \cE}$; and, \emph{ii)} the region of the injected powers is given by an affine transformation of $\cF$. Specifically, the first property \emph{i)} follows from the fact that power flows on different lines are decoupled; that is, it is possible to modify the angle $\theta_{mn}^{\phi}$ of a line $(m,n)$, while preserving the angle difference $\theta_{kl}^{\phi}$ of any other line $(k,l) \neq (m,n)$.    

Suppose now that the off-diagonal elements of $\bZ_{mn}$ are not zero; that is, $\bZ_{mn}^{-1} \neq (g_{mn} - j b_{mn})\bI_{|\cP_{nm}|}$. The total power flowing from node $m$ to node $n$ is given by $P_{m\rightarrow n} = \Re\{\bv_m^\cH\bZ_{mn}^{-1}(\bv_m - \bv_n)\}$, and it is now a function of $\{\theta_{mn}^{\phi}\}_{\phi \in \cP_{mn}}$, as well as of the angle differences $\{\theta_{m}^{\phi,\rho} := \theta_{m}^{\phi} - \theta_{m}^{\rho}\}_{\phi,\rho \in \cP_m}$ and $\{\theta_{n}^{\phi,\rho} := \theta_{n}^{\phi} - \theta_{n}^{\rho}\}_{\phi,\rho \in \cP_n}$. Different from the balanced case, the power flows $\{P_{m\rightarrow n}\}$ are no longer decoupled across lines. In fact, it is impossible to adjust the angles $\{\theta_{mn}^{\phi},\theta_{m}^{\phi,\rho},\theta_{n}^{\phi,\rho}\}$ to obtain a new flow on line $(m,n)$, without affecting the angle differences $\{\theta_{lm}^{\phi},\theta_{l}^{\phi,\rho},\theta_{m}^{\phi,\rho}\}$ and $\{\theta_{nk}^{\phi},\theta_{n}^{\phi,\rho},\theta_{k}^{\phi,\rho}\}$ for one of the other lines $(l,m)$ and $(n,k)$ connected to the nodes $m$ and $n$, respectively. Thus, the results of~\cite{LavaeiGeom} and~\cite{Tse12} no longer apply in the unbalanced setup.  

An analytical characterization of the flow region $\cF$ in the unbalanced case is challenging because of the number of voltage angles involved and the aforementioned coupling of the line power flows. Nevertheless, the following simple examples illustrate why one should expect a rank $1$ solution from the relaxed OPF even in an unbalanced setup.

%%%%%%%%%%%%%%%%%%%%%
\begin{figure}
	\centering
  \subfigure[Nonconvex problem.]{\includegraphics[width=0.23\textwidth]{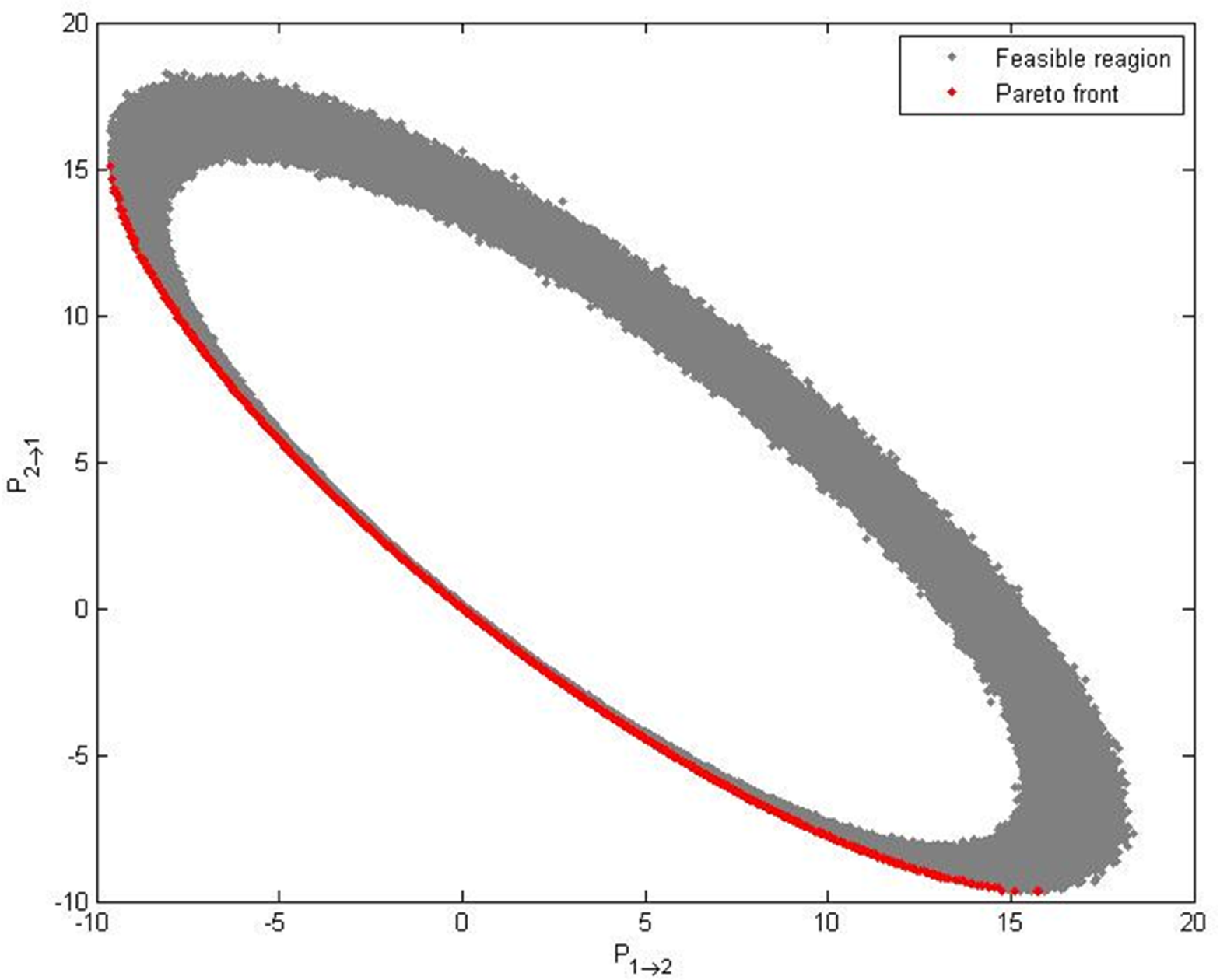}}
  \subfigure[Relaxed problem.]{\includegraphics[width=0.23\textwidth]{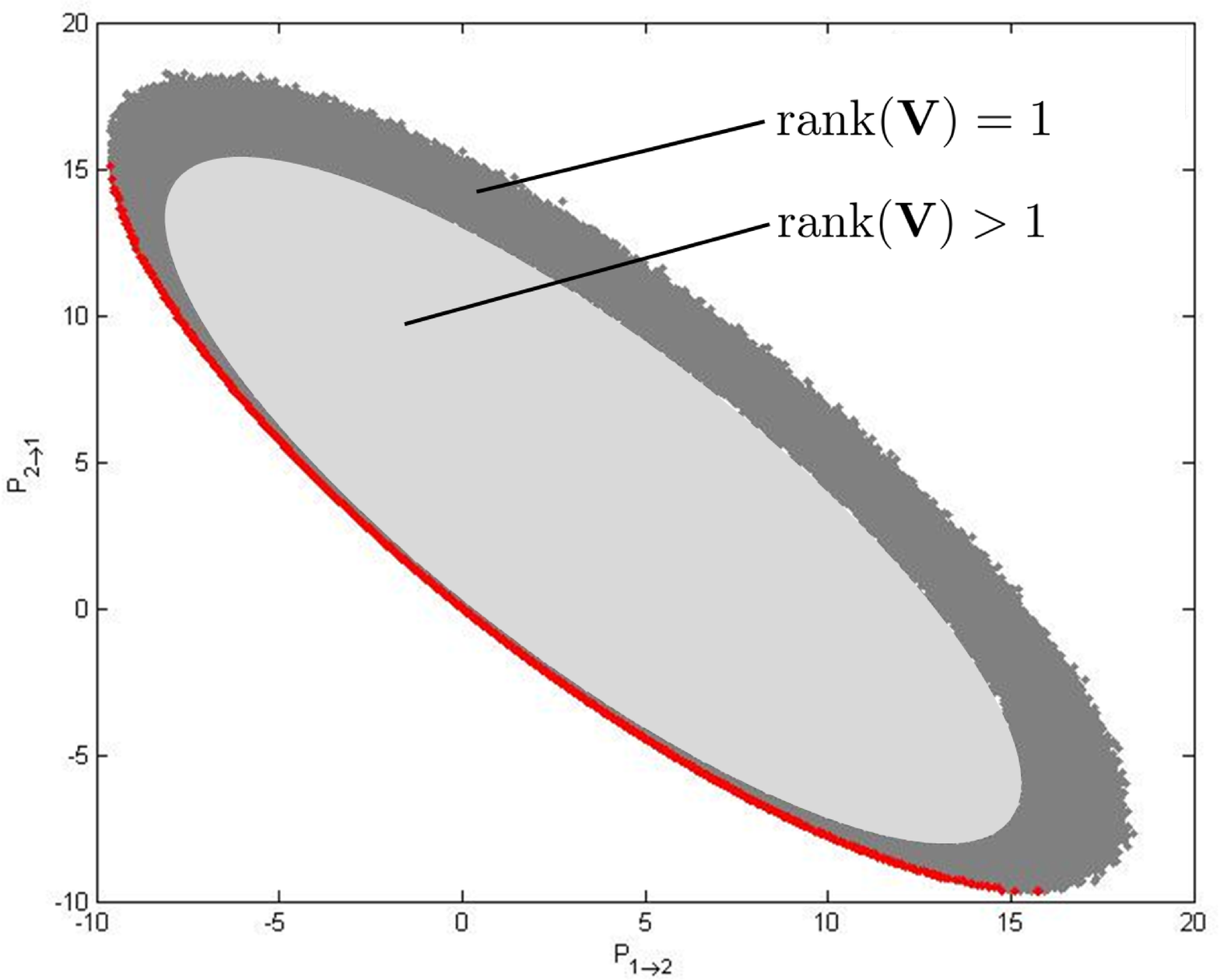}}
  \caption{Line flow region and its Pareto optimal points.}
  \label{fig:line}
\end{figure}
%%%%%%%%%%%%%%%%%%%%
%%%%%%%%%%%%%%%%%%%%%
\begin{figure}
	\centering
  \subfigure[Power injected at 3 nodes.]{\includegraphics[width=0.23\textwidth]{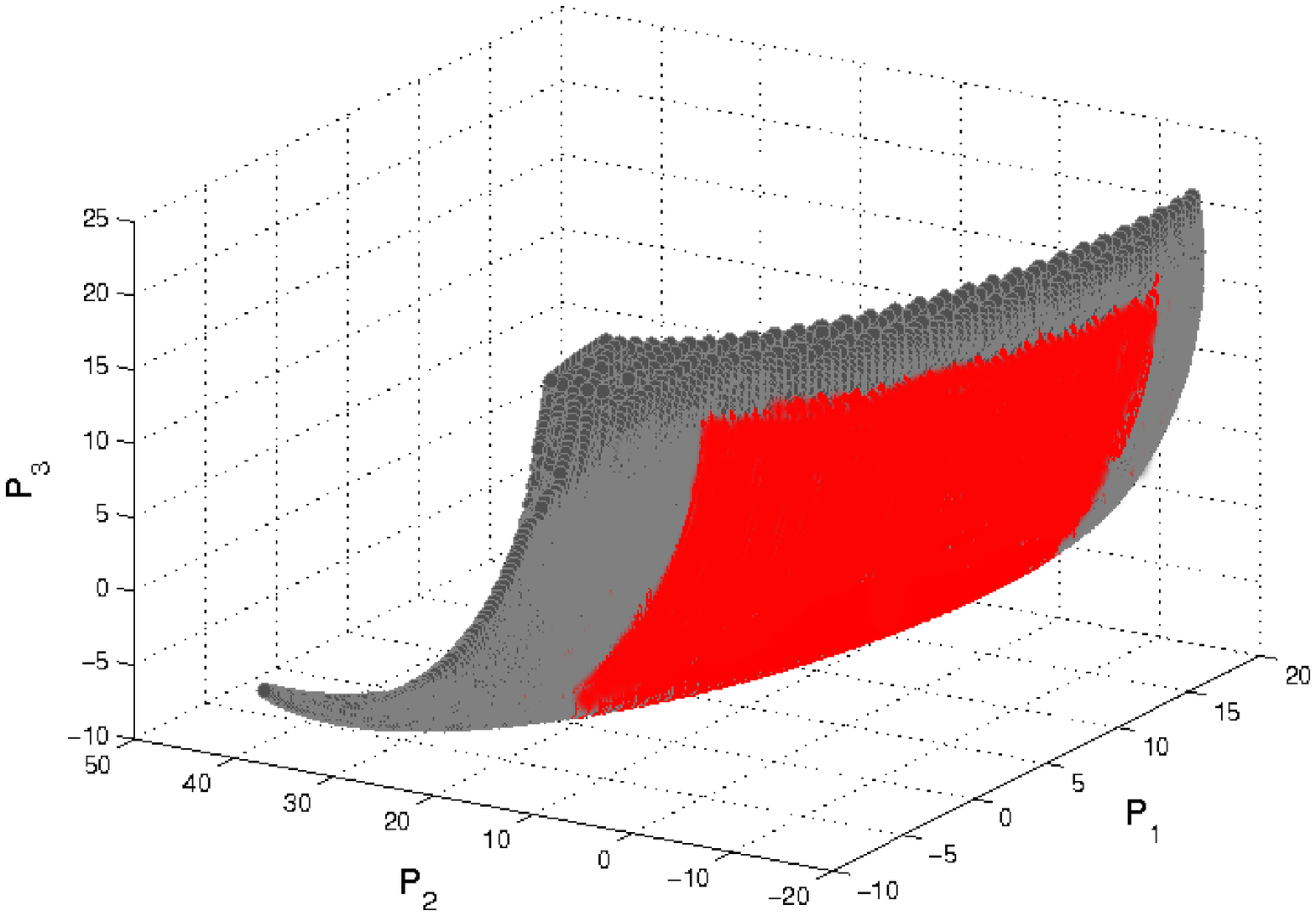}}
  \subfigure[Injection at 2 and 3 ($P_1 = -5$).]{\includegraphics[width=0.23\textwidth]{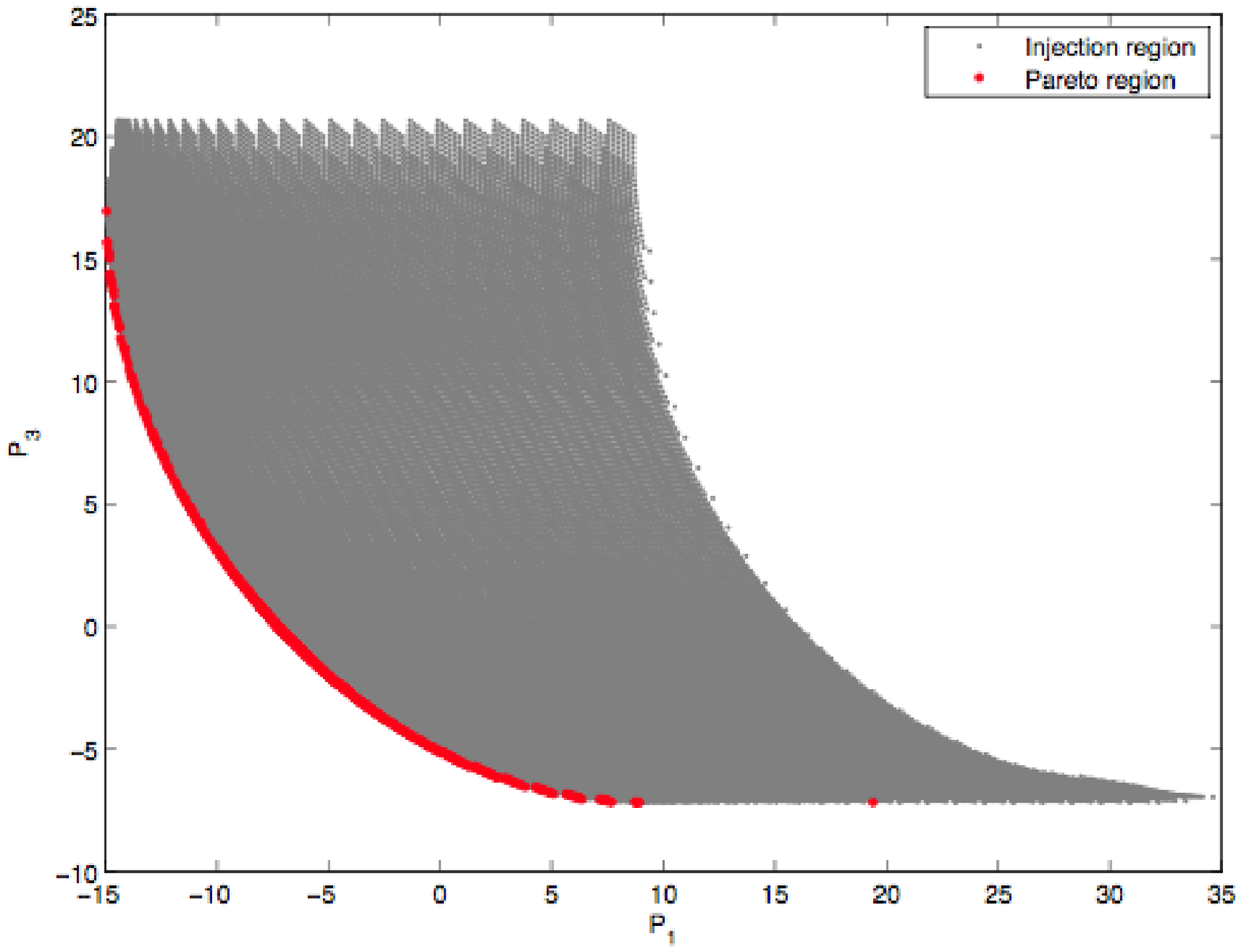}}
  \caption{Feasible power injection region and its Pareto optimal points.}
  \label{fig:injectedpowers}
\end{figure}
%%%%%%%%%%%%%%%%%%%%

Consider a 2-node unbalanced network, and suppose that a two-phase line connects the two nodes.
Let $\bZ_{12} = [(0.0753+j0.1181), (0.0156+ j 0.0502); (0.0156+ j 0.0502), (0.0744+ j 0.1211)]$ (as in~\cite{testfeeder}) and $\bY_{mn}^{(s)} = \mathbf{0}$, and assume that $|V_m^\phi| = 1$ for $n = 1,2$ and $\phi = a,b$. 
Fig.~\ref{fig:line}(a) depicts the feasible region $\cF_{12}$ of powers $(P_{1\rightarrow 2},P_{2\rightarrow 1})$ for $\theta_1^a = 0^\circ$, $\theta_{12}^{a} \in [-180^\circ,180^\circ]$, $\theta_{1}^{ab} \in [110^\circ,130^\circ]$, and $\theta_{2}^{ab} \in [110^\circ,130^\circ]$. Notice that $\cF_{12}$ is given by the Minkowski sum of the per-phase regions $ \{(P_{1\rightarrow 2}^{\phi},P_{2\rightarrow 1}^{\phi})\}$.   
It can be seen that the feasible flow region $\cF_{12}$ (the dark gray area) is a perturbed ellipsoid; specifically, expanding the expression of $P_{m\rightarrow n}$, it follows that $\theta_{1}^{ab}$ and $\theta_{2}^{ab}$ entail a perturbation of the center and of the axes of the ellipsoid that would be obtained if the line was balanced. Next, consider minimizing  a strictly increasing function  $C(P_{1\rightarrow 2},P_{1\rightarrow 2})$ over $\cF_{12}$. It follows that the set of Pareto optimal points are the ones represented by the red dots in Fig.~\ref{fig:line}(a). Consider now the convex hull of $\cF_{12}$, which amounts to augmenting the dark gray area with the light gray one shown in Fig.~\ref{fig:line}(b). Clearly, points belonging to the light gray area lead to a solution of (P3) with rank higher than 1. Notably, the Pareto points of $\cF_{12}$ and the ones of the relaxed region in Fig.~\ref{fig:line}(b) coincide; since $C(\cdot)$ is strictly increasing, the solution of (P3) must be on the Pareto boundary; thus, the optimal solution of (P3) has rank $1$, and it is an optimal solution also for the nonconvex problem (P2). Granted that the Pareto points in Fig.s~\ref{fig:line}(a) and (b) are the same,~\cite[Lemma~5]{LavaeiGeom} can be used to show that the two Pareto regions coincide also when constraints on the voltages are involved.  

The region of injected powers for a 3-node network is examined next. Let $\bZ_{12} = \bZ_{23} = [(0.0753+j0.1181), (0.0156+ j 0.0502); (0.0156+ j 0.0502), (0.0744+ j 0.1211)]$, and assume that $|V_n^\phi| = 1$ for all $n = 1,2,3$ and $\phi = a,b$.
Finally, let $P_{n} := P_{n}^a + P_{n}^b$, with $P_{n}^\phi$ the power injected at node $n = 1,2,3$ and phase $\phi = a,b$. The gray area depicted in Fig.~\ref{fig:injectedpowers}(a) corresponds to the region of feasible power injections $\cI$ for $-\tan^{-1}(\Re\{[\bZ_{12}^{-1}]_{1,1}\}/\Im\{[\bZ_{12}^{-1}]_{1,1}\}) < \theta_{nm}^{a} < \tan^{-1}(\Re\{[\bZ_{12}^{-1}]_{1,1}\}/\Im\{[\bZ_{12}^{-1}]_{1,1}\})$~\cite{LavaeiGeom}, and for the line-line angles
$\theta_{1}^{ab}, \theta_{2}^{ab}, \theta_{3}^{ab} \in [110^\circ,130^\circ]$; specifically, $\tan^{-1}(\Re\{[\bZ_{12}^{-1}]_{1,1}\}/\Im\{[\bZ_{12}^{-1}]_{1,1}\}) \approx 58^\circ$. The angle difference $\theta_{3}^{ab}$ is confined in the set $[110^\circ,130^\circ]$, as higher values are not likely to happen in practice~\cite{testfeeder}. Clearly, the gray region is nonconvex. If one considers minimizing a strictly increasing function of the injected powers $\{P_1, P_2,P_3\}$ over $\cI$, it follows that the Pareto optimal points are the ones color coded red in Fig.~\ref{fig:injectedpowers}(a). It can be noticed that the Pareto region does not change if one takes the convex hull of $\cI$; therefore, the solution of (P3) has rank 1, and it is an optimal solution for (P2). Next, fix the power injected (or absorbed) at node $3$, and consider minimizing a strictly increasing function of $P_1$ and $P_2$. The two-dimensional region of feasible powers $P_1$ and $P_2$ is depicted in gray in Fig~\ref{fig:injectedpowers}(b). Again, the Pareto front of the gray region and the one of its convex hull coincide.             
   
These examples suggest that the nonconvex OPF problem (P2) and its relaxed counterpart (P3) share their optimal solution when $-\tan^{-1}(\Re\{[\bZ_{mn}^{-1}]_{\phi,\phi}\}/\Im\{[\bZ_{mn}^{-1}]_{\phi,\phi}\}) < \theta_{nm}^{a} < \tan^{-1}(\Re\{[\bZ_{mn}^{-1}]_{\phi,\phi}\}/\Im\{[\bZ_{mn}^{-1}]_{\phi,\phi}\})$, and the angle differences between conductors $\theta_m^{\phi \rho}$ are small enough. This further motivates efforts toward analytical characterization of the power injection region in unbalanced distribution systems.

%%%%%%%%%%%%%%%%%%%%%%%%%%%%%%%%%%%%%%%%%%%%%%%
\section{Distributed Solution}
\label{sec:distributed}
%%%%%%%%%%%%%%%%%%%%%%%%%%%%%%%%%%%%%%%%%%%%%%%
 
Albeit polynomial, the computational complexity incurred by standard interior-point solvers for SDP does not scale well with the number of nodes $N$, and the number of constraints $N_c$~\cite{luospmag10}. However, this lack of scalability is incurred also by alternative methods based on off-the-shelf solvers for nonlinear programs~\cite{PaudyalyISGT}. Furthermore, the communication overhead required to collect data from all end users at the MEM, and subsequently disseminate the OPF solution back to the LCs, may lead to traffic congestions and substantial delays in the data delivery. Therefore, a \textit{distributed} SDP solver, with minimal computational and communication costs, is well motivated. A distributed approach is also desirable in order to address possible concerns regarding data privacy and integrity, and 
when the microgrid includes single- or multi-facilities that are managed independently from the rest of the network in order to pursue specific economic interests~\cite{Guerrero12}. 

Consider partitioning the microgrid into $L$ areas $\{\cA^{(l)} \subset \cN\}_{l = 1}^{L}$,  and suppose that each area is controlled by an LC. In the distributed SDP solver to be derived, each LC will solve an OPF problem of reduced dimension for its controlled area. Let $\bar{\cA}^{(l)}$ be an ``extended'' area defined as $\bar{\cA}^{(l)} := \cA^{(l)} \cup \{n | (m,n) \in \cE, m \in \cA^{(l)}, n \in \cA^{(i)}, l \neq j \}$; that is, $\bar{\cA}^{(l)}$ collects also the nodes belonging to different areas that are connected to $\cA^{(l)}$ by a distribution line. Sets $\{\bar{\cA}^{(l)}\}$ can be interpreted as counterparts of the ``regions'' considered in~\cite{Baldick97} in the context of transmission systems. Based on the overlaps among $\{\bar{\cA}^{(l)}\}_{l = 1}^{L}$, define the set of neighboring areas for the $l$-th one as $\bar{\cN}^{(l)} := \{j| \bar{\cA}^{(l)} \cap \bar{\cA}^{(j)} \neq 0\}$. Finally, 
let the vector $\bar{\bv}^{(l)}$ stack the complex line-to-ground voltages of the nodes 
in $\bar{\cA}^{(l)}$ (that is, $\{\bv_n\}_{n \in \bar{\cA}^{(l)}}$), and let $\bPhi_{P,n}^{\phi,l},\bPhi_{Q,n}^{\phi,l},\bPhi_{V,n}^{\phi,l},\bPhi_{m \rightarrow n}^{l}$ and $\bV^{(l)} := \bar{\bv}^{(l)} (\bar{\bv}^{(l)})^\cH$ denote the sub-matrices of
$\bPhi_{P,n}^{\phi},\bPhi_{Q,n}^{\phi},\bPhi_{V,n}^{\phi},\bPhi_{m \rightarrow n}$ and $\bV$, respectively, formed by extracting rows and columns corresponding to nodes
in $\bar{\cA}^{(l)}$.  With these notational conventions, it is possible to re-write the SDP (P3) as     
\begin{subequations}
\label{SDPdecoup}
\begin{align} 
& \hspace{-2.7cm} \mathrm{(P4)} \hspace{1.2cm} \min_{\bV} \,\, \sum_{l = 1}^{L} \tilde{C}_m^{(l)}(\bV^{(l)}) \label{mg4-cost}   \\
\textrm{s.t.} \,\, & \bV^{(l)} \in \cB^{(l)} \, ,  \,\, l = 1, \ldots L\label{mg4-region}  \\
& \bV \succeq \mathbf{0} \label{mg4-semipos}  
\end{align}
\end{subequations}
where $\cB^{(l)}$ denotes the set of sub-matrices $\bV^{(l)}$ satisfying the following constraints per area $\cA^{(l)}$:\begin{subequations}
\label{feasiblearea}
\begin{align} 
& \trace(\bPhi_{P,n}^{\phi,l} \bV^{(l)}) + P_{L,n}^{\phi} = 0, \hspace{.7cm} \forall\,\, \phi, \, \forall \, n \in  \cA^{(l)} \cap \cS  \label{mg4-P} \\
& \trace(\bPhi_{Q,n}^{\phi,l} \bV^{(l)}) + Q_{L,n}^{\phi} - y_{C,n}^{\phi}\trace(\bPhi_{V,n}^{\phi,l} \bV^{(l)})  = 0, \nonumber \\
& \hspace{3.8cm} \forall\,\, \phi, \, \forall \, n \in \cA^{(l)} \backslash (\cA^{(l)} \cap \cS ) \label{mg4-Q} \\
&P_{G,s}^{\textrm{min}} \leq \trace(\bPhi_{P,s}^{\phi,l} \bV^{(l)}) + P_{L,s}^{\phi} \leq P_{G,s}^{\textrm{max}},  \forall\,\, \phi, \, \forall \, n \in \cS^{(l)}  \label{mg4-Pg} \\
&Q_{G,s}^{\textrm{min}} \leq \trace(\bPhi_{Q,s}^{\phi,l} \bV^{(l)}) + Q_{L,s}^{\phi} \leq Q_{G,s}^{\textrm{max}},   \forall\,\, \phi, \, \forall \, n \in \cS  \label{mg2-Pg} \\
& (V_n^{\mathrm{min}})^2 \leq \trace(\bPhi_{V,n}^{\phi,l} \bV^{(l)}) \leq (V_n^{\mathrm{max}})^2, \forall\,\, \phi, \, \forall \, n \in \cA^{(l)} \label{mg2-voltage} 
\end{align}
\end{subequations}
with the additional constraint $ [\bV^{(l)}]_{\cP_0,\cP_0} = \bv_0 \bv_0^{\cH}$ for the area that contains the PCC. Further, 
$\tilde{C}_m^{(l)}(\bV^{(l)})$ represents the cost associated with the area $\cA^{(l)}$; that is, it collects the terms in~\eqref{eq:refcost1} (if $m = 1$) or~\eqref{eq:refcost1} (if $m = 2$) that pertain to $\cA^{(l)}$. For example, if one wishes to minimize the distribution losses, the corresponding cost becomes 
\begin{align} 
\tilde{C}_1^{(l)}(\bV^{(l)}) = & \sum_{(m,n) \in \cE | m,n \in \cA^{(l)}} \trace(\bPhi_{m \leftrightarrow n}^{l} \bV^{(l)}) \nonumber \\
& \hspace{-.9cm} + \sum_{(m,n) \in \cE | m \in \cA^{(l)}, n \in \bar{\cA}^{(l)} \backslash \cA^{(l)}} \frac{1}{2}  \trace(\bPhi_{m \leftrightarrow n}^{l} \bV^{(l)})
\end{align}
with $\bPhi_{m \leftrightarrow n}^{l} := \bPhi_{m \rightarrow n}^{l} + \bPhi_{n \rightarrow m}^{l}$. On the other hand, the expression for $\tilde{C}_2^{(l)}(\bV^{(l)})$ is not unique, depending on the specific network topology and operational setup. For instance, if ones assumes that each area is formed by one lateral or one sub-lateral~\cite{Guerrero12}, then $\tilde{C}_2^{(l)}(\bV^{(l)})$ accounts for the cost of power flowing into the (sub-)lateral from the ``backbone'' of the microgrid, and the cost of power generated within the (sub-)lateral. Instead, if customers are allowed to use DG units only to satisfy their own needs, $\tilde{C}_2^{(l)}(\bV^{(l)})$ boils down to $\tilde{C}_2^{(l)}(\bV^{(l)}) = c_0 \tilde{C}_1^{(l)}(\bV^{(l)}) + \sum_{s \in \cA^{(l)} \cap \cS} c_s \sum_{\phi \in \cP_s} P_{G,s}^{\phi}$. Either way, the equivalent formulation~\eqref{SDPdecoup} 
effectively expresses the cost as the superposition of local costs, and divides  network constraints on a per-area basis. However, even with such a decomposition the main challenge lies in the PSD constraint~\eqref{mg4-semipos}  that couples local matrices $\{\bV^{(l)}\}$. Indeed, if all submatrices $\{\bV^{(l)}\}$ were non-overlapping, the PSD constraint on $\bV$ would simplify to $\bV^{(l)} \succeq \mathbf{0}$ per area $l$, and (P4) would be decomposable in $L$ sub-problems. However, this equivalence fails to hold here, since submatrices share entries of $\bV$. The idea is to identify valid network topologies (that is, valid partitions of the microgrid in smaller areas) for which the PSD constraint decomposition is feasible. 

To this end, define the following two auxiliary graphs: 

\emph{i)} a ``macro'' graph $\cG_{\cA}$, where nodes represent the areas and edges are defined by the neighborhood sets $\{\bar{\cN}^{(l)}\}_{l = 1}^{L}$; and, 

\emph{ii)} a ``micro'' graph $\cG_{\cN}$ induced by the sub-matrices $\{\bV^{(l)}\}_{l = 1}^{L}$; that is, a graph with $\sum_{n \in \cN} |\cP_n|$ nodes (one per phase and node), with an edge connecting the nodes representing the voltages $V_n^{\phi}$ and $V_m^{\theta}$ if the entry of $\bV$ corresponding to $V_n^{\phi}(V_m^{\theta})^*$ is contained in one of the sub-matrices $\{\bV^{(l)}\}_{l = 1}^{L}$.

\noindent Examples of macro graphs are provided in Figs.~\ref{fig:F_feeder37} and~\ref{fig:F_toy_network}. 

Based on these auxiliary graphs, 
results on completing partial Hermitian
matrices will be leveraged to obtain PSD ones~\cite{Grone84}. These results rely on the so-termed chordal property of the graph $\cG_{\cN}$ induced by $\{\bV^{(l)}\}_{l = 1}^{L}$ to establish the equivalence between
positive semidefiniteness of the overall matrix $\bV$ and that of all
submatrices corresponding to the graph's maximal cliques.
%A similar technique has been used recently in \cite{Lam11} to distribute 
%the OPF problem in balanced transmission networks. 
Towards decomposing the PSD constraint into local ones, the following assumptions are made, which naturally suggest valid  partitions of the microgrid: 
\begin{enumerate}
\item
[\emph{(As1)}] The graph $\cG_{\cA}$ is a tree; and, 
\item
[\emph{(As2)}] $|\bar{\cA}^{(l)} \backslash (\bar{\cA}^{(l)} \bigcap \bar{\cA}^{(i)} ) | > 0$ for all $i,l = 1, \ldots, L$, $i \neq l$; that is, no nested extended areas are present.
\end{enumerate}
Condition (As1) is quite reasonable in tree distribution networks; for example, an area can be 
formed by a pair of nodes that are connected by a distribution line~\cite{Tse12}, or by laterals and sub-laterals~\cite{Guerrero12}. (As2) is a technical condition ensuring that the subgraph induced by $\bV^{(l)}$ is a maximal clique of $\cG_{\cN}$, which allows using the results of~\cite{Grone84}. Based on these assumptions, the following can be readily proved.

\vspace{.1cm}

\begin{proposition}
\label{prop:chordal}
Under (As1) and (As2), the graph $\cG_{\cN}$ is chordal, meaning that each of its cycles comprising four or more nodes has a chord. Furthermore,  
all its maximal cliques correspond to the elements of $\{\bV^{(l)}\}_{l = 1}^{L}$. 
\hfill $\Box$
\end{proposition}

\vspace{.1cm}

As established in~\cite{Grone84}, the PSD matrix $\bV$ is ``completable'' if
and only if $\cG_{\cN}$ is chordal, and all its submatrices 
corresponding to the maximal cliques of $\cG_{\cN}$ are PSD. Therefore, constraining 
$\bV$ to be PSD is tantamount to enforcing
the constraint on all local matrices $\bV^{(i)} \succeq \mathbf{0}$, $\forall l = 1,\ldots, L$. 
Notice that (As1) requires the macro graph $\cG_{\cA}$ to be a tree, while no conditions are 
imposed on the topology of the microgrid. Thus, it may be possible to find network partitions with an associated  chordal graph $\mathcal{G}_{\mathcal{N}}$ also in the case of weakly-meshed microgrids. 

Next, let $\cP_{lj} := \bigcup_{n \in (\bar{\cA}^{(l)} \cap \bar{\cA}^{(j)})} \cP_n$ collect the indexes corresponding to the voltages $\{\{V_n^{\phi}\}_{\phi \in \cP_n}\}$ of the nodes that the extended areas $\bar{\cA}^{(l)}$ and $\bar{\cA}^{(j)}$ share. For example, if areas $1$ and $2$ share nodes $n = 5$ and $n = 6$, then $\cP_{12}$ indexes the voltages $\{V_5^{\phi}\}_{\phi \in \cP_5}$ and $\{V_6^{\phi}\}_{\phi \in \cP_6}$. Further, define as $\bV^{(l)}_j$ the submatrix of $\bV^{(l)}$ collecting the rows and columns of $\bV^{(l)}$ corresponding to the voltages in $\cP_{l,j}$. With these definitions, and assuming that (As1) and (As2) hold, problem (P4) can be re-written in the following equivalent form:     
\begin{subequations}
\label{SDPdecoup}
\begin{align} 
& \hspace{-1.1cm} \mathrm{(P5)} \hspace{.5cm} \min_{\{\bV^{(l)}\}} \,\, \sum_{l = 1}^{L} \tilde{C}_m^{(l)}(\bV^{(l)}) \label{mg5-cost}   \\
\textrm{s.t.} \,\, & \bV^{(l)} \in \cB^{(l)} \, , \hspace{.60cm}  l = 1, \ldots L\label{mg5-region}  \\
&  \bV^{(l)}_j = \bV^{(j)}_l, \hspace{.5cm}  j \in \bar{\cN}^{(l)}, \,\, l = 1, \ldots, L  \label{consensus} \\
& \bV^{(l)} \succeq \mathbf{0} \, , \hspace{.95cm} l = 1, \ldots L  \label{mg5-semipos}  
\end{align}
\end{subequations}
where constraint~\eqref{consensus} enforces neighboring areas to consent on the entries of $\bV^{(l)}$ and $\bV^{(j)}$ that they have in common.  Clearly, constraints \eqref{consensus} couple the optimization problems across areas. To enable a fully distributed solution, consider introducing the auxiliary variables $\{\bW^{(l)}_j\}_{j \in  \bar{\cN}^{(l)}}$ and $\{\bX^{(l)}_j\}_{j \in  \bar{\cN}^{(l)}}$ per  area. With these auxiliary variables, (P5) can be equivalently re-stated as
\begin{subequations}
\label{SDPdecoup2}
\begin{align} 
& \hspace{-1.8cm} \mathrm{(P6)} \hspace{.7cm} \min_{\{\bV^{(l)} \succeq \mathbf{0} \}} \sum_{l = 1}^{L} \tilde{C}_m^{(l)}(\bV^{(l)}) \label{mg6-cost}   \\
\textrm{s.t.} \,\,  \bV^{(l)} & \in \cB^{(l)} \, ,   \hspace{.9cm} l = 1, \ldots L\label{mg6-region}  \\
 \Re\{\bV^{(l)}_j\} & = \bW^{(j)}_l, \hspace{.6cm}  j \in \bar{\cN}^{(l)}, \,\, l = 1, \ldots, L  \label{mg6-real} \\
  \Im\{\bV^{(l)}_j\} & = \bX^{(j)}_l, \hspace{.7cm}  j \in \bar{\cN}^{(l)}, \,\, l = 1, \ldots, L  \label{mg6-imag} \\
 \bW^{(j)}_l & = \bW^{(l)}_j , \,\,  \hspace{.55cm} j \in \bar{\cN}^{(l)}, \,\, l = 1, \ldots, L  \\
 \bX^{(j)}_l & = \bX^{(l)}_j , \,\,  \hspace{.65cm}  j \in \bar{\cN}^{(l)}, \,\, l = 1, \ldots, L . 
 \end{align}
\end{subequations}
 
A similar approach was followed by
\cite{Lam11}, which utilized either primal or dual iterations to
distribute the OPF in balanced transmission networks. A distributed OPF for  
balanced distribution feeders was derived in~\cite{Tse12}, where 
the dual (sub-)gradient ascent was used.  
Unfortunately, sub-gradient ascent methods do not always lead to a satisfactory solution; 
when the dual function is non-differentiable and the step size is fixed,
dual and primal iterates converge only on the average. What is more, 
recovering the primal variables from the optimal dual variables is not always guaranteed~\cite[Sec. 5.5.5]{BoVa04}. 
Besides dealing with unbalanced power networks, the novelty here consists in solving (P6) distributedly 
by resorting to the ADMM~\cite[Sec.~3.4]{BeT89}, a powerful scheme that has been successfully applied to distributed optimization and estimation in several contexts.
To this end, let $\{\bGamma^{(l)}_j\}$ and $\{\bLambda^{(l)}_j\}$ be the multipliers associated 
with~\eqref{mg6-real} and~\eqref{mg6-imag}, respectively,  and 
consider the partial quadratically-augmented Lagrangian of~\eqref{SDPdecoup2} as
\vspace{-.5cm}

{\small
\begin{align}
& \hspace{-.2cm} \cL(\{\bV^{(l)}\},\{\bW^{(l)}_j\},\{\bX^{(l)}_j\},\{\bGamma^{(l)}_j\},\{\bLambda^{(l)}_j\})  = \sum_{l = 1}^{L} \left\{ \tilde{C}_m^{(l)}(\bV^{(l)})  \right.
 \nonumber \\
& \hspace{0cm} + \sum_{j \in \bar{\cN}^{(l)}} \Big[ \trace \left( (\bGamma^{(l)}_j)^{\cT}(\Re\{\bV^{(l)}_j\}  - \bW^{(j)}_l)   \right)  
  \nonumber \\
& \hspace{.5cm} \left. + \trace \left( (\bLambda^{(l)}_j)^{\cT}(\Im\{\bV^{(l)}_j\}  - \bX^{(j)}_l)   \right) \right. \nonumber \\
& \hspace{.5cm} \left.  + \frac{\kappa}{2}\|\Re\{\bV^{(l)}_j\}  - \bW^{(j)}_l\|_F^2 + \frac{\kappa}{2}\|\Im\{\bV^{(l)}_j\}  - \bX^{(j)}_l\|_F^2  \Big] \right\} 
\label{eq:Lagrangian}
\end{align}}
\normalsize

\noindent where $\kappa \in \mathbb{R}^+$ is a positive constant~\cite[Sec.~3.4]{BeT89}. Then, the 
ADMM amounts to iteratively performing the following steps ($i$ denotes the iteration index): 

\vspace{.1cm}

\noindent \textbf{[S1] Update primal variables}:
\vspace{-.4cm}

\begin{align}
& \{\bV^{(l)}(i+1)\} := \nonumber  \\
& \arg \min_{\{\bV^{(l)}\succeq \mathbf{0}}  \cL(\{\bV^{(l)}\},\{\bW^{(l)}_j(i)\},\{\bX^{(l)}_j(i)\}, \nonumber \\
& \hspace{4cm} \{\bGamma^{(l)}_j(i)\},\{\bLambda^{(l)}_j(i)\}) \nonumber \\
& \hspace{1cm} \mathrm{s.t.} \,\, \bV^{(l)} \in \cB^{(l)}, \,\,  l = 1, \ldots, L. 
\label{eq:S1}
\end{align}

\noindent \textbf{[S2] Update auxiliary variables}:
\vspace{-.4cm}

\begin{align}
& \{\bW^{(l)}_j(i+1),\bX^{(l)}_j(i+1)\} := \nonumber  \\
& \arg \min_{\{\bW^{(l)}_j,\bX^{(l)}_j\}}  \cL(\{\bV^{(l)}(i+1)\},\{\bW^{(l)}_j\},\{\bX^{(l)}_j\}, \nonumber \\ 
& \hspace{4.3cm} \{\bGamma^{(l)}_j(i)\},\{\bLambda^{(l)}_j(i)\}) \nonumber \\
& \hspace{1cm} \mathrm{s.t.} \,\,  \bW^{(j)}_l = \bW^{(l)}_j ,  \,\, \bX^{(j)}_l  = \bX^{(l)}_j , \,\,  j \in \bar{\cN}^{(l)}, \,\, \,l
\label{eq:S2}
\end{align}

\noindent \textbf{[S3] Update dual variables}:
\vspace{-.4cm}

{\small
\begin{align}
\bGamma^{(l)}_j (i+1) & = \bGamma^{(l)}_j (i)  + \kappa (\Re\{\bV^{(l)}(i+1)\}  - \bW^{(j)}_l (i+1)) \label{eq:S3a} \\
\bLambda^{(l)}_j (i+1) & = \bLambda^{(l)}_j (i)  + \kappa (\Im\{\bV^{(l)}(i+1)\}  - \bX^{(j)}_l (i+1)) \label{eq:S3b}
\end{align}
}
\normalsize
In step [S1], the per-area matrices $\{\bV^{(l)}(i)\}$ are obtained by minimizing~\eqref{eq:Lagrangian}, where variables $\{\bW^{(l)}_j(i+1),\bX^{(l)}_j(i+1)\}$ and the multipliers are kept fixed to their previous iteration values. Likewise, the auxiliary variables are updated in [S2] by fixing $\{\bV^{(l)}(i+1)\}$ to their up-to-date values. Finally, the dual variables are updated in [S3] via dual sub-gradient ascent.

Interestingly, the ADMM iterations [S1]--[S3] can be simplified by exploiting the favorable decomposability of the Lagrangian. To this end, the following lemma is first needed. 

\vspace{.2cm}

\begin{lemma}
\label{lemma:dualupdates}
If the multipliers are initialized as $\bGamma^{(l)}_j (0) = \bLambda^{(l)}_j (0)= \mathbf{0}_{|\cP_{lj}|}$, then for every pair of neighboring areas $l$ and $j$ it holds that $\bGamma^{(l)}_j (i) + \bGamma^{(j)}_l (i) = \mathbf{0}_{|\cP_{lj}|}$ for each $i \geq 1$. Likewise, $\bLambda^{(l)}_j (i) + \bLambda^{(j)}_l (i) = \mathbf{0}_{|\cP_{lj}|}$ for each $i \geq 1$ if $\bLambda^{(l)}_j (0) = \bLambda^{(l)}_j (0)= \mathbf{0}_{|\cP_{lj}|}$. 
\end{lemma}
\emph{Proof.} See the Appendix. \hfill $\Box$

\vspace{.2cm}

Using Lemma~\ref{lemma:dualupdates}, steps [S1]--[S3] can be simplified as follows. Furthermore, convergence to the solution of the centralized SDP (P3) is established.    

\vspace{.2cm}

\begin{proposition}
\label{prop:ADMM}
If $\bGamma^{(l)}_j (0) = \bLambda^{(l)}_j (0)= \mathbf{0}_{|\cP_{lj}|}$, then [S1]--[S3] boil down to the following primal-dual updates:

\vspace{.1cm}

\noindent \textbf{[S1$^\prime$] Update} $\bV^{(l)}(i+1)$ \textbf{per area} $l = 1, \ldots, L$ as:
\begin{subequations}
\label{Pmg-per-area}
\begin{align} 
 \hspace{-.8cm} (\mathrm{P7}^{(l)}) & \hspace{.2cm} \min_{\bV^{(l)} \succeq \mathbf{0}, \{\alpha_{j} \geq 0,\beta_{j} \geq 0\}} \,\,  \cL^{(l)} (\bV^{(l)}, \{\alpha_{j},\beta_{j}\}) \label{cost-per-area} \\
& \textrm{s.t.} \,\, \bV^{(l)} \in \cB^{(l)} \\
& \hspace{.5cm} 
\left[ \begin{array}{cc}
-\alpha_{j} & \br_{j,\Re}^\cT \\
\br_{j} & -\bI
\end{array} \right] \preceq \mathbf{0}, \,\,\, \forall j \in \bar{\cN}^{(l)} \label{penalty} \\
& \hspace{.5cm} 
\left[ \begin{array}{cc}
-\beta_{j} & \br_{j,\Im} ^\cT \\
\br_{j,\Im}  & -\bI
\end{array} \right] \preceq \mathbf{0}, \,\,\, \forall j \in \bar{\cN}^{(l)} \label{penalty} 
\end{align}
\end{subequations}
where the local Lagrangian~\eqref{cost-per-area} is given by
\vspace{-.4cm}

{\small
\begin{align} 
& \cL^{(l)}(\bV^{(l)}, \{\alpha_{j},\beta_{j}\}) :=  \tilde{C}_m^{(l)}(\bV^{(l)}) + \sum_{j \in \bar{\cN}^{(l)}} \frac{\kappa}{2}(\alpha_j + \beta_j) \nonumber \\ 
& + \sum_{j \in \bar{\cN}^{(l)}} \left[ \trace \left( (\bGamma^{(l)}_j(i))^{\cT}\Re\{\bV^{(l)}_j\} \right)+ \trace \left( (\bLambda^{(l)}_j(i))^{\cT}\Im\{\bV^{(l)}_j\} \right) \right] \nonumber 
\end{align}}
\normalsize

\noindent and the vectors $\br_{j,\Re}$ and $\br_{j,\Im}$ collect the real and imaginary parts, respectively, of the entries of the matrix $\bV^{(l)}_j - \frac{1}{2}\left(\bV^{(l)}_j(i) + \bV^{(j)}_l(i) \right)$.  

\vspace{.2cm}

\noindent \textbf{[S2$^\prime$] Update dual variables locally per area} $l = 1,\ldots,L$:

\vspace{-.3cm}

{\small
\begin{align}
\hspace{-.3cm} \bGamma^{(l)}_j (i+1) & = \bGamma^{(l)}_j (i) + \frac{\kappa}{2} \left( \Re\{\bV^{(l)}_j(i+1)\} - \Re\{\bV^{(j)}_l(i+1) \right) \label{dualupdate1} \\
\hspace{-.3cm} \bLambda^{(l)}_j (i+1) & = \bLambda^{(l)}_j (i) + \frac{\kappa}{2} \left( \Im\{\bV^{(l)}_j(i+1)\} - \Im\{\bV^{(j)}_l(i+1) \right) \label{dualupdate2} 
\end{align}
}
\normalsize

\noindent Furthermore, for any $\kappa > 0$ the iterates $\{\bV^{(l)}(i)\}$, $\{\bGamma^{(l)}_j (i),\bLambda^{(l)}_j (i)\}$  produced by [S1$^\prime$]--[S2$^\prime$] are convergent, and $\lim_{i \rightarrow + \infty} \bV^{(l)}(i) = \bV^{(l)}_{\mathrm{opt}}$ for all $l = 1,\ldots, L$, with $\{\bV^{(l)}_{\mathrm{opt}}\}$ sub-matrices of the optimal solution $\bV_{\mathrm{opt}}$ of (P3).    
\end{proposition}
\emph{Proof.} See the Appendix. \hfill $\Box$

\vspace{.2cm}

At each iteration, the LC of area $l$ receives from the LCs of its neighboring areas $j \in \bar{\cN}^{(l)}$ matrices $\bV^{(j)}_l(i)$, and updates the local multipliers $\{\bGamma^{(l)}_j (i),\bLambda^{(l)}_j (i)\}$ via~\eqref{dualupdate1}--\eqref{dualupdate2}. These multipliers are \emph{locally} stored at area $l$, and they are not exchanged among LCs (in contrast, multipliers are exchanged per iteration in~\cite{Lam11}). Then, LC $l$ updates $\bV^{(l)}(i+1)$ by solving (P7$^{(l)}$), and transmits $\bV^{(l)}_j(i+1)$ to its neighboring areas $j \in \bar{\cN}^{(l)}$. % The overall distributed procedure is tabulated as Algorithm 1.

% 
%%%%%%%%%%%%%%%%%%%%%%%%%%%%%%%%%%%%%%%%
%\begin{algorithm}[t]
%\label{alg:decentralized}
%\caption{ADMM-based distributed solver} \small{
%\begin{algorithmic}
%
%\STATE Set $\bV^{(l)}$, and $\bGamma^{(l)}_j(0) = \bGamma^{(l)}_j(0) = \mathbf{0}$ for all $j \in \bar{\cN_0^\ell}$.
%
%\FOR {$i = 1,2,\ldots$ (repeat until convergence)} 
%
%\STATE 1: Update $\bV^{(l)(i)}$ by solving $(P7^{(l)})$. 
%
%\STATE 2: Send $\bV^{(l)_j(i)}$ to neighboring area $j$ Repeat for all $j \in \bar{\cN}^{(l)}$.
%
%\STATE 3: Receive $\bV^{(j)_l(i)}$ from neighboring areas $j \in \bar{\cN}^{(l)}$. 
%
%\STATE 4: Update multipliers via~\eqref{dualupdate1}--\eqref{dualupdate2}. 
%
%\ENDFOR
%\end{algorithmic}}
%\end{algorithm}
%%%%%%%%%%%%%%%%%%%%%%%%%%%%%%%%%%%%%%%

%%%%%%%%%%%%%%%%%%%%%%%%%%%%%%%%%%%%%%%%%%%%%%%
\section{Numerical Tests}
\label{sec:results}
%%%%%%%%%%%%%%%%%%%%%%%%%%%%%%%%%%%%%%%%%%%%%%%

%%%%%%%%%%%%%%%%%%%%
\begin{figure}[t]
\begin{center}
\includegraphics[width=0.4\textwidth]{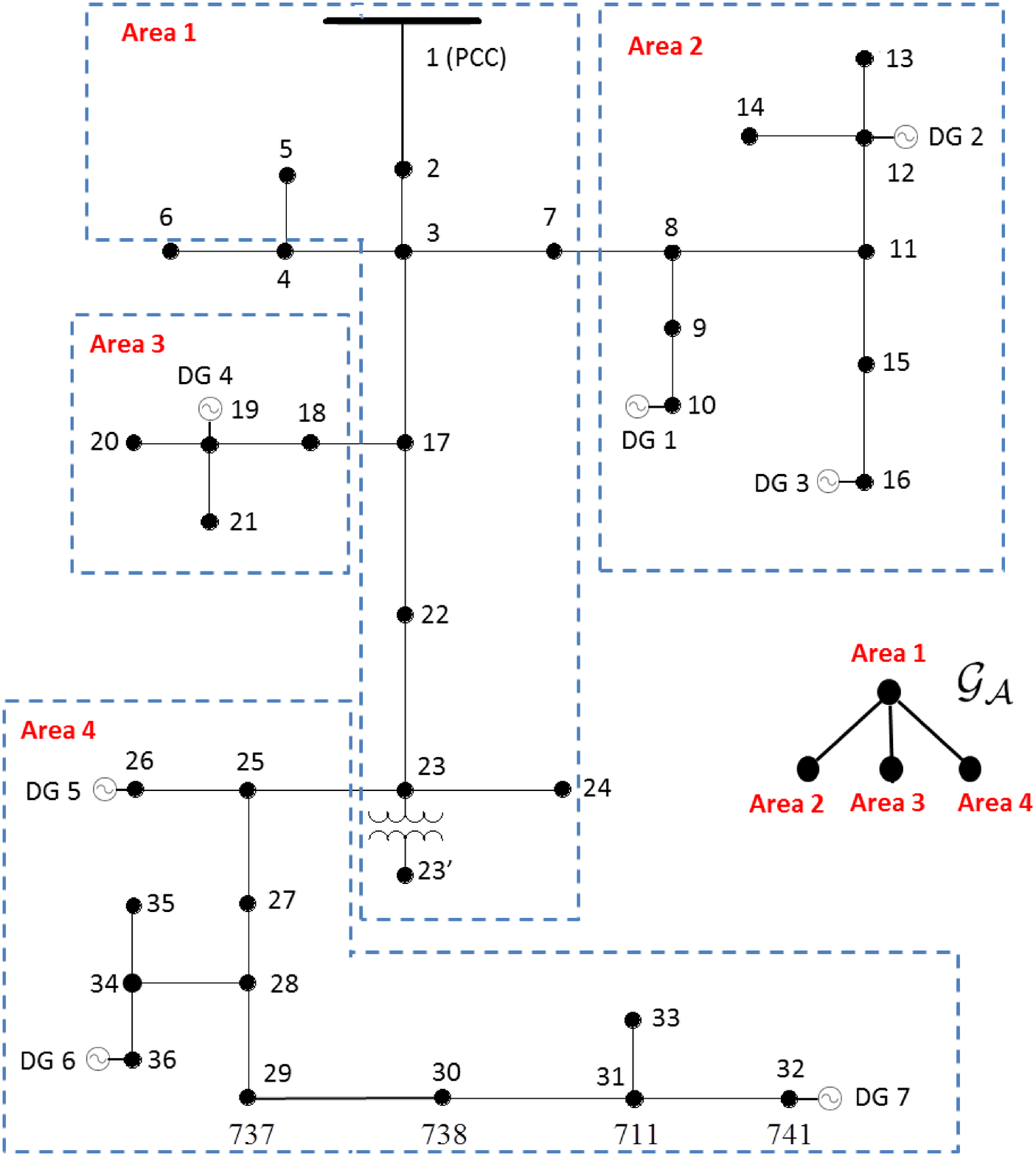}
\caption{First Test: IEEE 37-node feeder. }
\label{fig:F_feeder37}
\end{center}
\end{figure}
%%%%%%%%%%%%%%%%%%%%
%%%%%%%%%%%%%%%%%%%%
\begin{figure}[t]
\begin{center}
\includegraphics[width=0.3\textwidth]{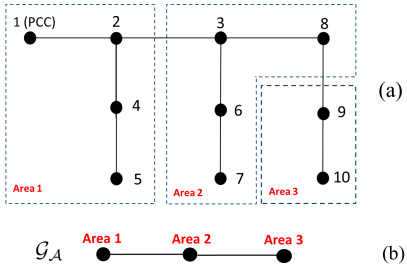}
\caption{Second test: (a) 10-node microgrid; (b) corresponding graph $\cG_{\cA}$.}
\label{fig:F_toy_network}
\end{center}
\vspace{-.5cm}
\end{figure}
%%%%%%%%%%%%%%%%%%%%
 
The SDP-based solver is tested here using the following two networks, operating in a grid-connected mode:

\begin{itemize}
	\item  the IEEE 37-node test feeder~\cite{testfeeder} shown in Fig.~\ref{fig:F_feeder37}; and
	\item  the 10-node 3-phase network depicted in Fig.~\ref{fig:F_toy_network}.
\end{itemize}
The optimization package \texttt{CVX}\footnote{[Online] Available: \texttt{http://cvxr.com/cvx/}}, along with the interior-point based solver \texttt{SeDuMi}~\cite{SeDuMi} are employed to implement the centralized and distributed solvers in \texttt{MATLAB}. 
 
The $4.8$ kV network of Fig.~\ref{fig:F_feeder37}  
is an actual portion of power a distribution network located in California; all the demanded complex powers are ``spot'' loads, and the network loading is very unbalanced~\cite{testfeeder}. Compared to the original scheme however, $7$ DG units are placed at nodes $\cS = \{10, 12, 16, 19, 26, 32, 36\}$. Specifically, single-phase conventional DG units supply a maximum real power of $50$ kW per phase, and are operated at a unit power factor (PF); that is, $Q_{G,s}^{\textrm{min}} = Q_{G,s}^{\textrm{max}} = 0$ for all $s \in \cS$. Line impedances and shunt admittances are computed based on the dataset in~\cite{testfeeder}. Finally, delta-wye conversions are performed whenever necessary. As for the network of Fig.~\ref{fig:F_toy_network}(a), the line admittances are all set to $\bZ_{mn} = [0.0693+j0.2036, 0.0312+j0.1003, 0.0316+j0.0847; 0.0312+j0.1003, 0.0675+j0.2096,0.0307+j0.0770;0.0316+j0.0847,0.0307+j0.0770,0.0683+j0.2070] \, \Omega$ (see~\cite{testfeeder}), which gives rise to an unbalanced operation of the network. Shunt admittances are neglected. Single-phase conventional DG units are placed at nodes $\cS = \{5,7\}$, and they can supply a maximum real power of $50$ kW per phase, at unit PF. All the loads are ``spot'', and the loading is assumed balanced. Specifically, the active and reactive loads are set to $\{0, 50, 0, 130, 130, 110, 110,   0,  90, 90,\}$ kW and $\{0, 20, 0,  82,  82,  60,  60,   0,  30, 30,\}$ kVAr, respectively, on each phase. The voltage magnitude at the PCC is $4.16$ kV. With these two choices, the performance of the proposed approach can be assessed for different network sizes, line characteristics, loading, and different network partitions shown next.    
   
The minimum and maximum utilization and service voltages are set to $V_n^{\textrm{min}} = 0.95$ pu and $V_n^{\textrm{max}} = 1.05$ pu for all nodes. Thus, voltage regulation is enforced without requiring changes in the voltage regulator taps (as in e.g.,~\cite{Paudyal11}). Further, the voltage at the PCC is set to $\bv_1 = [1 \angle 0^\circ, 1 \angle -120^\circ, 1 \angle 120^\circ]^\cT$ pu. The average computational time required by \texttt{SeDuMi} to solve the centralized problem (P3) was $9.0$ sec and $0.3$ sec (machine with Intel Core i7-2600 CPU @ 3.40GHz), which is significantly lower than the time required by commercial solvers for non-linear programs (see e.g.,~\cite{PaudyalyISGT}).     
   
%%%%%%%%%%%%%%%%%%%%
\begin{table}[h]
\caption{Results for the IEEE 37-node feeder}
\begin{center}
\begin{tabular}{|c|ccccc}
$c_s/c_0$ & $P_{\textrm{loss}}$ [kW] & $C_2$ [$\$ $] & $P_{0}$ [MW] & $P_{G}$ [MW] & $\rank(\bV_{\textrm{opt}})$ \\
\hline\\
0 &  39.30 &  57.9  & 1.4463 & 1.0500 & 1\\
0.25 & 39.30 &  68.4 & 1.4463 & 1.0500 & 1\\
0.50 & 39.30  & 78.9 & 1.4463 & 1.0500 & 1\\
0.75 & 38.97  & 89.3 & 1.4469 & 1.0491 & 1 \\ 
1 & 36.60 & 99.7 & 1.5979 & 0.8957 & 1\\
1.25 & 57.22 & 102.6 & 2.3088 & 0.2054 & 1
\end{tabular}
\end{center}
\label{powers1}
\vspace{-.4cm}
\end{table}
%%%%%%%%%%%%%%%%%%%%

%%%%%%%%%%%%%%%%%%%%
\begin{table}[h]
\caption{Results for the 10-node network in Fig.~\ref{fig:F_toy_network}}
\begin{center}
\begin{tabular}{|c|ccccc}
$c_s/c_0$ & $P_{\textrm{loss}}$ [kW] & $C_2$ [$\$ $] & $P_{0}$ [MW] & $P_{G}$ [MW] & $\rank(\bV_{\textrm{opt}})$ \\
\hline\\
0 &   18.38 &  59.9  & 1.4984 &  0.3000 & 1\\
0.25 & 18.38 &  62.9 & 1.4984 &  0.3000 & 1\\
0.50 & 18.38  & 65.9 & 1.4984 &  0.3000 & 1\\
0.75 & 18.38  & 68.9 & 1.4984 &  0.3000 & 1 \\ 
1 & 18.27 & 71.9 & 1.5457 & 0.2526 & 1\\
1.25 & 23.08 & 72.1 & 1.8031 & 0.0000 & 1
\end{tabular}
\end{center}
\label{powers2}
\vspace{-.4cm}
\end{table}
%%%%%%%%%%%%%%%%%%%%

%%%%%%%%%%%%%%%%%%%%
\begin{figure}[t]
\begin{center}
\includegraphics[width=0.50\textwidth]{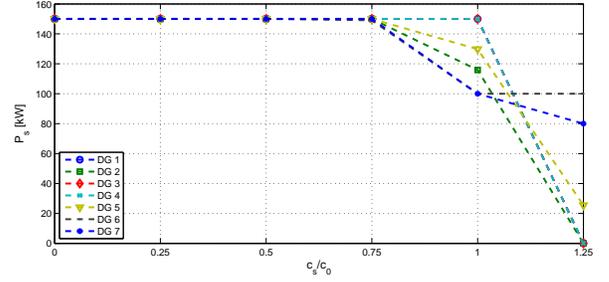}
\caption{Power generated by the DG units [kW].}
\label{fig:F_DG_Pgen}
\end{center}
\vspace{-.7cm}
\end{figure}
%%%%%%%%%%%%%%%%%%%%  

Tables~\ref{powers1} and~\ref{powers2} list the real power drawn at the PCC $P_0 := \sum_{\phi \in \cP_{0}} V_{0}^\phi (I_{0}^\phi)^*$, the total power generated by the DG units $P_G := \sum_s \sum_\phi P_{G,s}^{\phi}$, and the overall power losses $P_{\textrm{loss}}$ and costs of supplied power, when the cost~\eqref{econ_dispatch_cost} is employed. The costs of supplied power are set to $c_0 = 40$ $\$$/MW and $c_1 = \ldots = c_7 \in \{0,10,20,30,40,50\}$ $\$$/MW. Notice that minimizing~\eqref{econ_dispatch_cost} is tantamount to minimizing the power loss~\eqref{power_loss_cost} when $c_0 = 1$ and $c_s = 1$ for all $s \in \cS$. In fact, it  can be clearly seen in Tables~\ref{powers1} and~\ref{powers2} that the power loss is minimized for this choice of $\{c_s\}$. Powers $P_{\textrm{loss}}$, $P_0$, and $P_G$ remain the same when $c_s \leq c_0$; however, when $c_s \geq c_0$ (which holds during the late-night hours) the DG units reduce the generated active powers; further, the power loss becomes significantly higher when DG units are not used at the maximum extent.
Fig.~\ref{fig:F_DG_Pgen} depicts the active powers generated by DG units as a function of $c_s/c_0$, for the IEEE 37-node test feeder. It can be noticed that the DG units electrically close to the PCC are not utilized when $c_s \geq c_0$; on the other hand, DG 6 and 7 still operate at more than $50 \%$ of their maximum capacity.

%%%%%%%%%%%%%%%%%%%%
\begin{figure}[t]
\begin{center}
\includegraphics[width=0.5\textwidth]{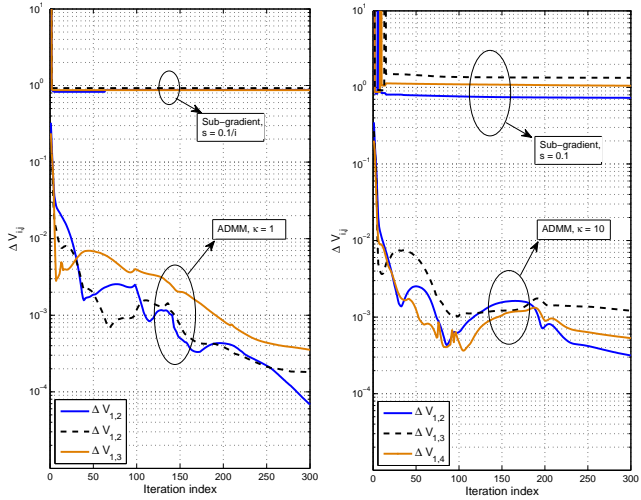}
\caption{Convergence of the ADMM, with feeder partitioned as in Fig.~\ref{fig:F_feeder37}.}
\label{fig:Convergence_feeder37_2}
\end{center}
\vspace{-.5cm}
\end{figure}
%%%%%%%%%%%%%%%%%%%%

%%%%%%%%%%%%%%%%%%%%
\begin{figure}[t]
\begin{center}
\includegraphics[width=0.5\textwidth]{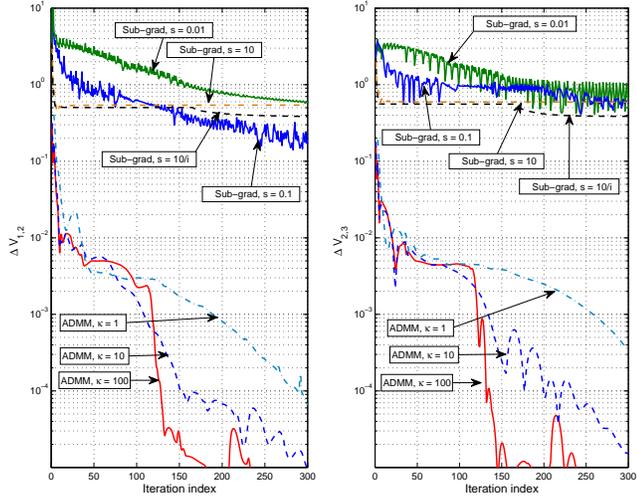}
\caption{Convergence of the ADMM, network partitioned as shown in Fig.~\ref{fig:F_toy_network}.}
\label{fig:Convergence_toy_network_2}
\end{center}
\vspace{-.5cm}
\end{figure}
%%%%%%%%%%%%%%%%%%%%

Interestingly, the rank of matrix $\bV_{\textrm{opt}}$ was \emph{always} 1. Therefore, the \emph{globally} optimal solutions of (P2) (and hence of the nonconvex (P1)) were always attained. In other words, no lower power losses or costs of supplied power can be attained with alternative OPF solution approaches. The rank of matrix $\bV_{\textrm{opt}}$ was greater than $1$ for the IEEE 37-node test feeder when $c_s/c_0 > 1.75$ (two non-zero eigenvalues, with $\lambda_1 \propto 10^2$ and $\lambda_2 \propto 10^{-2}$ ). Nevertheless, a rank-1 solution was readily obtained upon raising the voltage magnitude at the PCC from $1$ to $1.02$ pu. This further prompts an analytical characterization of the feasibility region of (P1). Additional tests were performed on the IEEE 13-node feeder~\cite{testfeeder}, and rank-1 matrices $\bV_{\textrm{opt}}$ were again always obtained (results are not reported here due to space limitation). 

Convergence of the proposed distributed SDP solver is showcased in Fig.~\ref{fig:Convergence_feeder37_2}, where the 37-node feeder is partitioned as shown in Fig.~\ref{fig:F_feeder37}. This partition resembles the case where laterals include multi-facilities that are managed independently from the rest of the network~\cite{Guerrero12}. The trajectories corresponding to  
$\Delta V_{l,j}(i) := \|\bV_l^{(j)}(i) - \bV_j^{(l)}(i)\|_1/36$ per iteration $i$, are reported for different values of the ADMM parameter $\kappa$, and are compared with the ones obtained by using the sub-gradient ascent-based distributed algorithm developed in~\cite{Tse12} ($s$ denotes the step size, which is assumed to be either constant, or monotonically decreasing). It can be noticed that the proposed distributed solver exibit a considerably faster convergence than the one based on the sub-gradient. Considering that the nominal voltage of the feeder is $4.8$ kV, in less then than $50$ iterations the average gap between the entries of $\bV_l^{(j)}(i)$ and $\bV_j^{(l)}(i)$ is on the order of a few volts.  Further, notice that the convergence rate is approximately the same for $\Delta V_{1,2}(i)$, $\Delta V_{1,3}(i)$, and $\Delta V_{1,4}(i)$.      
The average computational time required by \texttt{SeDuMi} to solve each sub-problem was of $6, 3, 0.3$, and $2.7$ sec. Lower   computational times can be obtained by selecting areas of smaller size.

Fig.~\ref{fig:Convergence_toy_network_2} illustrates the convergence of the distributed SDP solved when the considered 10-node network is partitioned as shown in Fig~\ref{fig:F_toy_network}. The ADMM-based method outperforms the one based on the sub-gradient here too. In this case, the gaps between the elements of matrices $\{\bV_l^{(j)}(i)\}$ rapidly vanish when $\kappa = 100$ and $\kappa = 10$ after approximately $100$ iterations.  In this case, the average time required by \texttt{SeDuMi} to solve each sub-problem was $0.1$ sec.

%%%%%%%%%%%%%%%%%%%%%%%%%%%%%%%%%%%%%%%%%%%%%%%
\section{Concluding remarks}
\label{sec:conclusions}
%%%%%%%%%%%%%%%%%%%%%%%%%%%%%%%%%%%%%%%%%%%%%%%

The OPF problem was considered for microgrids operating in an unbalanced setup. 
Inspite of the inherent non-convexity, the SDP relaxation technique was advocated
to obtain a convex problem. As corroborated by
numerical tests, the main contribution of the proposed approach
consists in offering the potential to obtain the globally optimal
solution of the original nonconvex OPF. A distributed SDP solver was also developed by resorting to the ADMM. The distributed algorithm ensures scalability with respect to the microgrid size, 
robustness to communication outages, and preserves data privacy and integrity.

%%%%%%%%%%%%%%%%%%%%%%%%%%%%%%%%%%%%%%%%%%%%%%
\appendix
%%%%%%%%%%%%%%%%%%%%%%%%%%%%%%%%%%%%%%%%%%%%%%

\emph{Proof of Lemma}~\ref{SDPreformulation}. 
Equality~\eqref{eq:V} can be readily established by noticing that $|V_n^{\phi}|^2 = \bv^\cH \bar{\be}_n^{\phi} (\bar{\be}_n^{\phi})^{\cT} \bv = \trace(\bar{\be}_n^{\phi} (\bar{\be}_n^{\phi})^{\cT}\bV)$. To prove~\eqref{eq:P}, notice first that the injected apparent power at node $n$ and phase $\phi$ is given by $V_n^{\phi} (I_n^{\phi})^* = (V_n^{\phi *} I_n^{\phi})^* = (\bv^\cH  \bar{\be}_n^{\phi} (\bar{\be}_n^{\phi})^{\cT} \bi)^{\cH}$. Next, using $\bi = \bY \bv$, it follows that $(\bv^\cH  \bar{\be}_n^{\phi} (\bar{\be}_n^{\phi})^{\cT} \bi)^{\cH} = (\bv^\cH \bar{\be}_n^{\phi} (\bar{\be}_n^{\phi})^{\cT} \bY \bv)^{\cH} = (\bv^\cH  \bY^{\phi}_n \bv)^{\cH} = \bv^\cH  (\bY^{\phi}_n)^{\cH} \bv$, which can be equivalently rewritten as $\trace(\bY_n^{\phi} \bV)$. Thus, the injected real and reactive powers can be obtained by using, respectively, the real and imaginary parts of $(\bY^{\phi}_n)^{\cH}$. Finally, to prove~\eqref{eq:Pmn}, define first the $|\cP_{mn}| \times 1$ vector $\bi_{m \rightarrow n}$ collecting the complex currents flowing from $m$ to $n$ on each phase, and notice that $\bi_{m \rightarrow n} = \bZ_{mn}^{-1}(\bv_m - \bv_n) = \bA_{m \rightarrow n} \bv$ and $\bv_m = \bB_m \bv$.   Then, it follows that $P_{m \rightarrow n} = \Re\{\bi_{m \rightarrow n}^\cH \bv\} = \Re\{\bv^\cH \bA_{m \rightarrow n}^\cH \bB_m \bv\} = \Re\{\trace(\bA_{m \rightarrow n}^\cH \bB_m \bV)\}$.

\vspace{.2cm}

\emph{Proof of Lemma}~\ref{lemma:dualupdates}
The proof is provided for the multipliers $\bGamma^{(l)}_j$; same steps can be followed 
for the dual variables $\bLambda^{(l)}_j$. 

Suppose that $i = 1$, and notice that the optimization problem~\eqref{eq:S2} decouples per pair of neighboring areas; specifically, a number of sub-problems is to be solved with respect to (wrt)  the only pair of variables $(\bW^{(j)}_l, \bW^{(l)}_j)$ for $j\in \bar{\cN}^{(l)}$ (and clearly $l\in \bar{\cN}^{(j)}$). Neglecting irrelevant terms,  and setting the first-order derivative of $\cL(\{\bV^{(l)}(1)\},\{\bW^{(l)}_j\},\{\bX^{(l)}_j\},\{\bGamma^{(l)}_j(0)\},\{\bLambda^{(l)}_j(0)\})$ wrt $(\bW^{(j)}_l, \bW^{(l)}_j)$ to zero, 
it readily follows that the minimizer $\bW^{(l)}_j(1)$ is given by
\begin{equation}  
\label{proofDual_1}
\bW^{(l)}_j(1) = \frac{1}{2}\left(\bV^{(l)}_j(1)  + \bV^{(j)}_l(1) \right).
\end{equation} 
Next, substituting~\eqref{proofDual_1} into [S3], and setting $\bGamma^{(l)}_j(0) = \bGamma^{(l)}_j(0) = \mathbf{0}$ yields 
$\bGamma^{(l)}_j(1)  = \frac{\kappa}{2} \left( \bV^{(l)}_j(1) - \bV^{(j)}_l(1) \right)$ and 
$\bGamma^{(j)}_l(1)  = \frac{\kappa}{2} \left( \bV^{(j)}_l(1) - \bV^{(l)}_j(1) \right)$.
It then follows by induction that $\bGamma^{(l)}_j(i) = - \bGamma^{(j)}_l(i)$ also for all subsequent iterations. 

\vspace{.2cm}

\emph{Proof of Proposition}~\ref{prop:ADMM}
Substituting~\eqref{proofDual_1} into [S3], and setting $\bGamma^{(l)}_j(0) = \bGamma^{(l)}_j(0) = \mathbf{0}$, one can readily obtain~\eqref{dualupdate1}--\eqref{dualupdate2}.  Consider now [S1$^\prime$]. First, substitute~\eqref{proofDual_1} into~\eqref{eq:Lagrangian}, and discard irrelevant terms. 
Next, introduce an auxiliary variable $\alpha_{j}$ to upper bound the squared $l_2$ norm of $\br_{j,\Re}$, and add the constraint $\|\br_{j,\Re}\|_2^2 \leq \alpha_{j}$. Thus,~\eqref{penalty} follows by Schur's complement. Indeed, it is certainly possible to introduce one single variable $\alpha$ to upper bound the squared $l_2$ norm of the vector $[\br_{1,\Re}^\cT,\br_{1,\Im}^\cT,\ldots,\br_{|\bar{\cN}^{(l)}|,\Re}^\cT,\br_{|\bar{\cN}^{(l)}|,\Im}^\cT]^\cT$, and derive a linear matrix inequality similar to~\eqref{penalty} by using Schur's complement.              
Finally, since the cost in (P6) is convex, and the constraint set is bounded and convex, every limit point of $\{\bV^{(l)}(i)\}$ is an optimal solution to (P6), as established in~\cite[Prop.~4.2]{BeT89}. Since (P6) and (P3) are equivalent under (As1) and (As2) hold, every limit point of $\{\bV^{(l)}(i)\}$ is also an optimal solution to (P3).

%%%%%%%%%%%%%%%%%%%%%%%%%%%%%%%%%%%%%%%%%%%%%%
\bibliographystyle{IEEEtran}
\bibliography{Biblio_pow_systems}
%%%%%%%%%%%%%%%%%%%%%%%%%%%%%%%%%%%%%%%%%%%%%%

\end{document}